\documentclass[11pt, reqno]{amsart}
\usepackage[active]{srcltx}
\usepackage {color}
\usepackage{amsmath,amsthm,amssymb}
\usepackage{epsfig}
\usepackage{graphicx}
\usepackage{amssymb}
\usepackage{latexsym}
\usepackage{pdfpages}
\usepackage{bm}
\usepackage{enumitem}

\usepackage{ulem} 

\usepackage{ifpdf}

\ifpdf 
  \usepackage[hidelinks]{hyperref}
\else 
\fi 

 \usepackage[usenames,dvipsnames]{pstricks}
 \usepackage{epsfig}
 \usepackage{pst-grad} 
 \usepackage{pst-plot} 

\usepackage{color}

\usepackage{subcaption}

\newtheorem{theorem}{Theorem}[section]
\newtheorem{lemma}[theorem]{Lemma}
\newtheorem{definition}[theorem]{Definition}
\newtheorem{proposition}[theorem]{Proposition}
\newtheorem{corollary}[theorem]{Corollary}
\newtheorem{remark}[theorem]{Remark}
\newtheorem{notation}[theorem]{Notation}

\newtheorem*{theorem*}{\it Theorem}

\DeclareMathOperator*{\argmin}{arg\,min}

\def\vint_#1{\mathchoice%
          {\mathop{\kern 0.2em\vrule width 0.6em height 0.69678ex depth -0.58065ex
                  \kern -0.8em \intop}\nolimits_{\kern -0.4em#1}}%
          {\mathop{\kern 0.1em\vrule width 0.5em height 0.69678ex depth -0.60387ex
                  \kern -0.6em \intop}\nolimits_{#1}}%
          {\mathop{\kern 0.1em\vrule width 0.5em height 0.69678ex
              depth -0.60387ex
                  \kern -0.6em \intop}\nolimits_{#1}}%
          {\mathop{\kern 0.1em\vrule width 0.5em height 0.69678ex depth -0.60387ex
                  \kern -0.6em \intop}\nolimits_{#1}}}
\def\vintslides_#1{\mathchoice%
          {\mathop{\kern 0.1em\vrule width 0.5em height 0.697ex depth -0.581ex
                  \kern -0.6em \intop}\nolimits_{\kern -0.4em#1}}%
          {\mathop{\kern 0.1em\vrule width 0.3em height 0.697ex depth -0.604ex
                  \kern -0.4em \intop}\nolimits_{#1}}%
          {\mathop{\kern 0.1em\vrule width 0.3em height 0.697ex depth -0.604ex
                  \kern -0.4em \intop}\nolimits_{#1}}%
          {\mathop{\kern 0.1em\vrule width 0.3em height 0.697ex depth -0.604ex
                  \kern -0.4em \intop}\nolimits_{#1}}}

\def\R{\mathbb R}
\def\N{\mathbb N}

\numberwithin{equation}{section}

\renewcommand{\theequation}{\arabic{section}.\arabic{equation}}
\def\1{\raisebox{2pt}{\rm{$\chi$}}}

\definecolor{violet(ryb)}{rgb}{0.53, 0.0, 0.69}
\usepackage{mathtools}
\mathtoolsset{showonlyrefs}

\newcounter{mgncount}

\begin{document}

\title[Anisotropic p-capacity]{\bf Characterization of the subdifferential and minimizers for the anisotropic p-capacity}

\author[E. Cabezas-Rivas, S. Moll, M. Solera]{  Esther Cabezas-Rivas, Salvador Moll and Marcos Solera}

\address{E. Cabezas-Rivas: Departament de Matem\`atiques,
Universitat de Val\`encia, Dr. Moliner 50, 46100 Burjassot, Spain.
 {\tt esther.cabezas-rivas@uv.es  }}

\address{S. Moll: Departament d'An\`{a}lisi Matem\`atica,
Universitat de Val\`encia, Dr. Moliner 50, 46100 Burjassot, Spain.
 {\tt j.salvador.moll@uv.es }}
\address{M. Solera: Departament d'An\`{a}lisi Matem\`atica,
Universitat de Val\`encia, Dr. Moliner 50, 46100 Burjassot, Spain. {\tt marcos.solera@uv.es } \vspace{0.2cm}\newline
\phantom{MMM. Solera:} Departamento de Matem\'aticas, Universidad Aut\'onoma de Madrid, C/
Francisco Tom\'as y Valiente, 7, Facultad de Ciencias, m\'odulo 17, 28049 Madrid, Spain.
 {\tt marcos.solera@uam.es  }
}


\date{\today. The work of the first author is partially supported by the AEI (Spain) and FEDER project PID2019-105019GB-C21, and by the GVA project AICO 2021 21/378.01/1. The second and third authors have been partially supported  by ``Conselleria d'Innovaci\'o, Universitats, Ci\`encia i Societat Digital'', project AICO/2021/223. The third author has also been supported by the Ministerio de Universidades (Spain) and NextGenerationUE, programme ``Recualificaci\'on del sistema universitario espa\~{n}ol'' (Margarita Salas), Ref. UP2021-044; the ``Conselleria d'Innovaci\'o, Universitats, Ci\`encia i Societat Digital'', programme ``Subvenciones para la contratación de personal investigador en fase postdoctoral'' (APOSTD 2022),  Ref. CIAPOS/2021/28; the MICINN (Spain) and FEDER, project PID2021-124195NB-C32; and by the ERC Advanced Grant 834728.}


\setcounter{tocdepth}{1}

%

\begin{abstract}
We obtain existence of minimizers for the $p$-capacity functional defined with respect to a centrally symmetric anisotropy for $1 < p<\infty$, including the case of a crystalline norm in $\R^N$. The result is obtained by a characterization of the corresponding subdifferential  and it applies for unbounded domains of the form $\R^N \setminus \overline{\Omega}$ under mild regularity assumptions (Lipschitz-continuous boundary) and no convexity requirements on the bounded domain $\Omega$. If we further assume an interior ball condition (where the Wulff shape plays the role of a ball), then any minimizer is shown to be Lipschitz continuous.

\end{abstract}


\maketitle


{ \renewcommand\contentsname{Contents }
\setcounter{tocdepth}{3}
 }

 \section{Introduction and statement of main results}

 \subsection{The p-capacity and its anisotropic version}

 For $N \geq 2$, the $p$-capacity of a given compact set $K\subset\R^N$ is defined as
 $${\rm Cap}_p(K)=\inf\left\{\int_{\R^N} \|\nabla u\|^p\,\mathrm{d}x \, :\, u\in C_c^\infty(\R^N), \, u\geq 1 {\rm \ in \ } K\right\},$$
 where $\|\cdot\|$ denotes the Euclidean norm and $C_c^\infty$ are smooth functions with compact support. The relevance of this quantity comes from its geometric meaning, since for $p>1$ the usual geometric functionals like area or volume do not provide a satisfactory depiction of properties of domains. Here is where inequalities involving the $p$-capacity naturally arise.

 From the viewpoint of partial differential equations, the search for minimizers of the $p$-capacity, called $p$-{\it capacitary functions} or equilibrium potentials, has attracted a lot of attention due to its relation to the $p$-Laplace equation. In fact, when $1<p<N$ and under suitable smoothness of the boundary of $K$, a unique $p$-capacitary function, $u\in L^{\frac{pN}{N-p}}(\R^N)$ and $\nabla u \in L^p(\R^N; \R^N)$, exists and it satisfies the Euler-Lagrange equation 
 $$\left\{\begin{array}{cl}-\Delta_p u=-{\rm div}(\|\nabla u\|^{p-2}\nabla u)=0 & {\rm in \ }\R^N\setminus K \smallskip \\
 \hspace*{-2.8cm} u=1 & {\rm in \ }  K \smallskip \\ \hspace*{-2.8cm} u \to 0 & {\rm as } \, \|x\| \to \infty. \end{array}\right.$$

 The purpose of this paper is to obtain existence of $p$-capacitary functions within anisotropic media, with regularity assumptions so mild as to include the case of a crystalline anisotropy. It is crucial to move out of classical Euclidean isotropic settings to allow the possibility of environments
 where properties differ with the direction, since this extra flexibility will be key for applications to crystal growth or noise removal.

 Accordingly, replacing the Euclidean norm in the definition of the $p$-capacity functional by a generic norm $F$ in $\R^N$, leads to the {\it anisotropic $p$-capacity functional}:
 $${\rm Cap}^F_p(K):=\inf\left\{\int_{\R^N} F^p(\nabla u)\,\mathrm{d}x \, : \, u\in C_c^\infty(\R^N), u\geq 1 {\rm \ in \ } K\right\}.$$

 In the range $1<p<N$, existence and uniqueness of minimizers (called anisotropic $p$-capacitary functions) have been shown under different conditions on the anisotropy $F$ (most of the results deal with $C^\infty(\R^N\setminus\{0\})$ uniformly elliptic norms) and the domain $K$. Moreover, in these cases the minimizer is the unique solution to the following PDE:
 \begin{equation} \label{smooth-anis-PDE}
 \left\{\begin{array}
 {cl} -{\rm div}(F^{p-1}(\nabla u)\nabla F(\nabla u))=0 & {\rm in \ }\R^N\setminus K \smallskip \\ u=1 & {\rm in \ }K \smallskip \\ u\to 0 & {\rm as \ }\|x\|\to \infty.
 \end{array}\right.
 \end{equation}

 Our main goal is to find a similar characterization of minimizers of the anisotropic $p$-capacity by relaxing the assumptions on regularity of $F$ and $K$ as much as possible, and allowing a non-necessarily constant boundary condition. To state our main result in this direction, we need to introduce some definitions and notation.

 \subsection{Characterization of minimizers by means of the subdifferential}

 Let $\Omega \subset \R^N$ be an open bounded set with Lipschitz-continuous bounda\-ry and $\varphi\in W^{1-\frac1{p}, p}(\partial\Omega)$. Set
 \[W^{1,p}_\varphi (\Omega) := \{u\in W^{1,p}(\Omega) \, : \, u=\varphi \text{ on } \partial \Omega\}.\]
 For a continuous norm $F$ in $\R^N$, we consider the energy functional
 $$\mathcal{F}_{\Omega,\varphi}:L^p(\Omega)\rightarrow [0,+\infty]$$  defined by
 \begin{equation} \label{defF}
 \mathcal{F}_{\Omega,\varphi}(u):=\left\{\begin{array}{ll}
 \displaystyle\int_{\Omega} F^p(\nabla u), & \hbox{if } u\in W^{1,p}_\varphi(\Omega),   \smallskip \\
 +\infty, & \hbox{otherwise}.
 \end{array}\right.
 \end{equation}
 Notice that, since the domain of the above functional is bounded and $F^p$ is convex and coercive, \cite[Theorem 5 in section 8.2.4]{Evans} ensures the existence of a minimizer. However, since $F$ (and therefore $F^p$) is not required to be strictly convex, uniqueness is not anymore guaranteed.

 Our objective is to characterize the minimizers of this energy functional by means of the corresponding Euler-Lagrange equation (notice that here we are not allowed to write $\nabla F$ as in \eqref{smooth-anis-PDE}). Since $\mathcal{F}_{\Omega,\varphi}$ is a proper and convex functional, we have that
 \begin{equation} \label{equiv}
 u\in \argmin\left\{ \mathcal{F}_{\Omega,\varphi}\right\} \quad \text{if, and only if,} \quad 0\in\partial\mathcal{F}_{\Omega,\varphi}(u),
 \end{equation}
 where the latter stands for the subdifferential of the energy functional, whose exact expression (and accordingly the concrete Euler-Lagrange equation for this problem) is unknown and  does not follow from standard arguments, due to the non smoothness of the integrand. Therefore, our first goal is to characterize the subdifferential; more precisely, we get
 \begin{theorem}[Characterization of the subdifferential in $\Omega$ bounded] \label{Emin}
 	Let \\ $1 < p < \infty$, $\Omega \subset \R^N$ a bounded domain with Lipschitz-continuous boundary and $\varphi \in W^{1-\frac1{p},p}(\partial \Omega)$. If $u\in W^{1,p}_\varphi(\Omega)$ and  $v\in L^{p'}(\Omega)$, then the following are equivalent:
 	\begin{itemize}
 		\item[(a)] $v  \in \partial\mathcal{F}_{\Omega,\varphi}(u)$.
 		\item[(b)] There exists $z\in L^\infty(\Omega;\R^N)$, with $F^\circ(z) \leq 1$ and $z \cdot \nabla u = F(\nabla u)$, such that $v=-{\rm div} (p F^{p-1}(\nabla u) z)$ in the weak sense, that is,
 		$$\int_{\Omega} vw = \int_{\Omega} p F^{p-1}(\nabla u) z\cdot\nabla w \quad \hbox{ for every } w\in W_0^{1,p}(\Omega).$$
 	\end{itemize}
 	Here $F^\circ$ denotes the dual of the norm $F$.
 \end{theorem}

 The strategy and methods carried out to prove \autoref{Emin}, suitably adapted to work with homogeneous Sobolev spaces, permit to obtain the corresponding characterization for unbounded exterior domains of the form $\R^N \setminus \overline \Omega$ (see Theorem \ref{Emin2}), where $\Omega\subset\R^N$ is again a bounded domain with Lipschitz boundary. As a byproduct of this result,
 we get a characterization for all minimizers of ${\rm Cap}_p^F(\overline\Omega)$ for bounded Lipschitz domains $\Omega$.

 \begin{corollary}
 	Let  $1 < p < \infty$, $\Omega \subset \R^N$ a bounded domain with Lipschitz boundary and $\varphi \in W^{1-\frac1{p},p}(\partial \Omega)$. Consider either $D = \Omega$ or $D = \R^N \setminus \overline \Omega$, then any minimizer $u$ of the energy functional $\mathcal{F}_{D,\varphi}$ is a weak solution of
 	\begin{equation}\label{Eq-in-out}
 	\left\{ \begin{array}{ll} {\rm div}\big( F^{p-1}(\nabla u) z\big) = 0 & \text{in }  D\smallskip \\
 	u = \varphi & \text{on }  \partial D, \smallskip \\
 	u \to 0 \text{ as } \|x\| \to \infty & \text{if }  D = \R^N \setminus \overline \Omega, \quad p < N\end{array} \right.
 	\end{equation}
 	for some $z\in L^\infty(D;\R^N)$, with $F^\circ(z) \leq 1$ and $z \cdot \nabla u = F(\nabla u)$. Here the equality on the boundary is understood in the sense of traces. In particular, if $\varphi \equiv 1$, then
 	$${\rm Cap}_p^F(\overline\Omega)=\inf\left\{\mathcal F_{\R^N\setminus\overline\Omega,1}(u) : u\in L^{\frac{pN}{N-p}}\left(\R^N\right)\right\},$$
 	and hence anisotropic $p$-capacitary functions are solutions to \eqref{Eq-in-out}.
 \end{corollary}

 Up to our knowledge, the most general result in this spirit can be found in \cite{BiCi} (see also \cite{Memo}). The authors show existence, uniqueness and regularity of a $p$-capacitary function $u\in C^{2,\alpha}(\R^N\setminus K)$ in the case that $K$ is convex and with boundary of class $C^{2,\alpha}$ for norms $F\in C^{2,\alpha}(\R^N\setminus\{0\})$ such that $F^p$ is twice continuously differentiable in $\R^N\setminus\{0\}$ with a strictly positive definite Hessian matrix.

Here we generalize the results  in the previous literature in several directions. Indeed, regarding the anisotropy, we only require it to be a norm, without any extra assumptions on smoothness or uniform ellipticity. In particular, we allow norms whose dual unit balls have corners and/or straight segments, including crystalline cases as the $\ell_\infty$ or $\ell_1$ norms in $\R^N$.

 Secondly, we study the Dirichlet problem both in bounded domains with Lipschitz boundary and in exterior domains, as well as for a generic Dirichlet boundary constraint. 
 Moreover, we do not require that the set $K$ is convex and, about its regularity, we just ask for Lipschitz continuity of the boundary, instead of the traditional $C^{2, \alpha}$ much stronger constraint. Lipschitz cannot be further weakened because it is the milder condition that guarantees e.g. that the trace operator is surjective on the fractional Sobolev space $W^{1-\frac1{p},p}$.

 In short, we will be working in an unfriendly setting in the sense that we cannot perform any argument that involves second derivatives of the functions $u$ nor principal curvatures (even in weak sense) of $\partial K$. To overcome these additional technical difficulties, the proof of \autoref{Emin} strongly relies on the theory of maximal monotone operators in Banach spaces as in the case of $p=1$, previously studied in \cite{M}.

 First, we associate a possibly multivalued operator $\mathcal A_\varphi: L^p(\Omega)\to L^{p'}(\Omega)$ to item $(b)$ in Theorem \ref{Emin}. In order to ensure that this operator coincides with $\partial\mathcal F_{\Omega,\varphi}$, we show that $\mathcal A_\varphi\subseteq \partial\mathcal F_{\Omega,\varphi}$ and that both are maximal monotone (see \autoref{AeqDF}). The maximal monotonicity of $\mathcal A_\varphi$ will be proved by verifying that the range condition $L^{p'}(\Omega)=R(\mathcal J_{L^p}+\mathcal A_\varphi)$ holds (cf. Proposition \ref{Range-0}), with $\mathcal J_{L^p}$ being the duality mapping.

 The latter, in turn, needs an approximation process with a sequence of coercive, monotone and weakly continuous operators defined on $W^{1,p}_\varphi(\Omega)$. For the continuity, we will approximate the norm $F$ by its Moreau-Yosida approximation while for the coercivity, we will add a $p$-Laplacian term. Suitable a-priori estimates and the use of the Minty-Browder technique (see Lemma \ref{Minty-tec}) will permit to pass to the limit, first in the Yosida regularization and then in the $p$-Laplacian to finally achieve that the range condition holds.

\subsection{Construction of minimizers with extra properties}

 	We also show that there exists one $p$-capacitary function which is trapped between two explicit solutions of the Euler-Lagrange equations. By translation invariance, we can suppose that $0\in\Omega$ and, since $\Omega$ is bounded, we can find $0<r_1<r_2$ such that $\mathcal W_{r_1}\subset\Omega\subset\mathcal W_{r_2}$, with $\mathcal W_r$ being the  Wulff shape with radius $r$; i.e., the ball of radius $r$ with respect to the dual norm of $F$. Then, there are three $p$-capacitary functions $u, u_{r_1}, u_{r_2}$, minimizers of ${\rm Cap}_p^F(\overline\Omega)$, ${\rm Cap}_p^F(\overline{\mathcal W_{r_1}})$ and ${\rm Cap}_p^F(\overline{\mathcal W_{r_2}})$, respectively, such that $u_{r_1}\leq u\leq u_{r_2}$. Moreover, $u_{r_1}$ and $u_{r_2}$ are explicitly given as follows:
 		\begin{theorem} \label{umin-trapped}
 		If $1<p<N$ and $\Omega \subset\R^N$ is a bounded domain with Lipschitz  boundary, there exists a minimizer $u$ of  ${\rm Cap}_p^F(\overline\Omega)$ such that
 		\begin{enumerate}
 			\item[{\rm (a)}] $u\in L^{p^\ast}(\R^N \setminus \overline \Omega)$ is a weak solution of \eqref{Eq-in-out} with $D = \R^N \setminus \overline \Omega$, $p^\ast:=\frac{Np}{N-p}$ and $\nabla u \in L^p(\R^N \setminus \overline \Omega; \R^N)$.
 			\item[{\rm (b)}] We can find constants $0 < r_1 < r_2$ with $\mathcal W_{r_1} \subset \Omega \subset  \mathcal W_{r_2}$ so that
 			\[u_{r_1} \leq u \leq u_{r_2} \qquad \text{with } \ u_{r_i} =(F^\circ)^{\frac{p-N}{p-1}}\left(\frac{\cdot}{r_i}\right), \quad i = 1,2.\]
 			\item[{\rm (c)}] If $F$ is required to be strictly convex, then the minimizer is unique.
 		\end{enumerate}
 	\end{theorem}
 	
 	In order to prove the above result, we need to approximate $p$-capacitary functions with minimizers of the {\it relative capacity} with respect to a ball. Recall that, in the anisotropic framework, Wulff shapes play the role of balls. Therefore, for $\overline\Omega\subset\mathcal W_{R}$, we consider $${\rm Cap}_p^F(\overline\Omega;\mathcal W_{R}):=\inf\left\{\int_{\mathcal W_R} F^p(\nabla u) : u\in C_0^\infty(\mathcal W_R)\,, u\geq 1 {\rm \ in \ }\Omega\right\},$$
 	which is also known as the anisotropic $p$-capacity of the condenser
 	$(\overline\Omega;\mathcal W_{R})$ or the condenser anisotropic $p$-capacity of the obstacle $\overline \Omega$  in the bounded domain $\mathcal W_R$. Notice that the unique minimizer clearly satisfies that $u|_{\Omega} \equiv 1$, so we actually work within the annular domain $\Omega_R := \mathcal W_R \setminus \overline \Omega$, and try to let $R \to \infty$, as in the classical isotropic setting \cite{Lewis}.
 	
 	If we assume that $F$ is strictly convex, a comparison argument directly yields two barriers $v_{r_1,R}$ (lower) and $v_{r_2,R}$ (upper) into which the minimizer $u_R$ to ${\rm Cap}_p^F(\overline\Omega;\mathcal W_{R})$ is trapped. Since $u_R$ are shown to be increasing with respect to $R$, we can pass to the limit when $R\to \infty$ and we get the result.
 	
 	In the case that $F$ is not a strictly convex norm,  we need to approximate it with a sequence of strictly convex norms in a uniform way; to obtain uniform bounds on the barriers and then finish the argument as in the strictly convex case. Notice that, under this generality, one does not expect to get uniqueness, as this does not happen in the extremal case $p=1$ (see Appendix \ref{nonU} for an example with infinitely many minimizers for $p=1$).
 	
 	\subsection{Regularity of minimizers}
 	We complete the paper with the study of the regularity of minimizers. We show that all minimizers both of the relative $p$-capacity with respect to a Wulff shape and of the anisotropic $p$-capacity are Lipschitz continuous, provided that the domain is regular enough in the following sense:
 	\begin{definition}[Uniform interior ball condition] \label{UIBC}
 		Let $r>0$. We say that $\Omega$ satisfies the $\mathcal W_r$-condition if, for any $x\in\partial \Omega$, there exists $y\in \R^N$ such that
 		\[\mathcal W_r+y\subseteq \overline{\Omega} \quad \text{ and } \quad x\in\partial \left(\mathcal W_r+y\right).\]
 	\end{definition}
 	\noindent This condition is milder than $F$-regularity for non-convex domains (see e.g. the discussion in Lemma 2.8 and Remark 2.9  of \cite{CCMN}). In this setting, we conclude
 	
 	\begin{theorem} \label{lip-thm}
 		Let $r>0$ and suppose that $\Omega$ satisfies the $\mathcal W_r$-condition. Then any minimizer of ${\rm Cap}_p^F(\overline\Omega;\mathcal W_{R})$  is Lipschitz continuous. Moreover,  any minimizer of the energy functional $\mathcal{F}_{\R^N \setminus{\overline \Omega},1}$  is also Lipschitz continuous.
 	\end{theorem}
 	
 \noindent We point out that this is the best expected regularity since the explicit solutions $u_{r}$ in Theorem \ref{umin-trapped} are only Lipschitz continuous in the case that $F$ does not have any extra regularity assumption. The proof requires an application of our comparison arguments, which is trickier than for the previous theorem, as we need that the upper and lower barrier coincide on the boundary, in order to exploit a regularity result from \cite{MT}.
 	
 	\subsection{Geometric and physical meaning of $p$-capacity} As it is the main character of this paper, let us talk briefly about the physical and geometric relevance of $p$-capacity. Physically speaking,  ${\rm Cap}_2(K)$ measures the total electric charge flowing into $\R^N \setminus K$ across the boundary $\partial K$. But this interpretation does not restrict to electric charges, it can also be applied to heat transfer or even fluid flow through a porous medium.
 	
 	Indeed, the problems studied above can be interpreted as the steady states of such flows. In the classical case ($p=2$), Ohm's law says that the electric current is driven by the field $J =- c \nabla u$, where $u$ is the corresponding $p$-capacitary potential, and $c$ denotes the conductivity. But all the physical laws (Ohm, Fourier or Darcy) governing the aforementioned flows are empirical and linearity is just a simplifying assumption, hence the next level of complexity should consider flows driven by $J = - c \|\nabla u\|^{p-2} \nabla u$ for $p > 1$, which has already been studied in the context of turbulent flows and deformation plasticity (pioneering works in this direction are \cite{Phi, AtK}).
 	
 	As suggested by Pólya in \cite{Pol}, the thermal analogy works as a source of geometric intuition:  take a body (e.g.~a cat) within a uniform infinite medium whose temperature vanishes at infinity, while the skin of the cat is kept at a constant temperature (that we normalize to 1). Then the thermal conductance (quantity of outgoing heat per time unit) is, up to a constant which depends on the nature of the ambient, equal to the electrostatic capacity of the cat. In addition, we have all noticed that, to protect themselves from the cold, cats tend to curl up in a ball; this happens in order to minimize  the thermal conductance or, equivalently, their capacity.
 	
 	This can be formalized by means of {\it isocapacitary inequalities} (see \cite{May}), telling that,
 	among all sets with fixed volume, balls  minimize $p$-capacity, i.e.,
 	\[{\rm Cap}_p(K) \geq {\rm Cap}_p(B_r)\]
 	where $r$ is such that the Lebesgue measure of $B_r$ coincides with that of $K$. Equality holds if, and only if, $K = B_r$, up to a set of zero $p$-capacity. There are further interesting characterizations of balls as the equality case of Minkowski type inequalities, which relate suitable powers of the $p$-capacity and integrals involving a $p$-power of the mean curvature of $\partial K$ (cf. \cite {Agos}).
 	
 	Anisotropic generalizations of the above interpretations come naturally by considering bodies embedded in non-uniform media. For the corresponding anisotropic inequalities, see \cite{Li} and the references therein.
 	

 	\subsection{Structure of the paper} The paper is organized as follows. We first introduce in section \ref{back} the basic background material about anisotropies, maximal monotone theory in Banach spaces and Yosida regularization in Hilbert ambients, while section
 	\ref{subdif-1} gathers all the approximation arguments needed to prove \autoref{Emin}. These include estimates for the Moreau-Yosida approximation of $F^p$ and the Yosida approximation
 	of $\partial F^p$ (Lemma \ref{aux-Fp}), as well as a result in the spirit of Minty-Browder (Lemma \ref{Minty-tec}) gi\-ving a sufficient condition for elements to belong to $\partial F^p$.   Then in
 	section \ref{unbdd-sec} we introduce the technical machinery of homogeneous Sobolev spaces to extend the characterization of the subdifferential to unbounded exterior domains (\autoref{Emin2}). Section \ref{comp-barr} includes a comparison result (Lemma \ref{ComparisonPrincipleForApproximations}) for strictly convex norms; the obtaining of minimizers $u_R$ for the relative $p$-capacity ${\rm Cap}_p^F(\overline\Omega;\mathcal W_{R})$ as limits of minimizers of energy functionals, which involve strictly convex anisotropies approximating a generic norm $F$ (Lemma \ref{seq-min}); and the proof that  $u_R$ are trapped between explicit solutions of the corresponding Euler-Lagrange problem within annular domains (Proposition \ref{comp-res}). Then the next natural step is to carry over these barrier arguments in bigger and bigger rings to reach \autoref{umin-trapped} as the outer boundary tends to infinity, which stands for the content of section \ref{sec-ult}; in turn,
 	 \autoref{lip-thm} is proven in section \ref{reg-sect}. Finally, we include an appendix with an explicit example to justify the lack of uniqueness in the case $p=1$.

 \section{Notation and background material} \label{back}

 \subsection{Anisotropies and Wulff shape} \label{Anisot}

A continuous function $F:\R^N \rightarrow [0, \infty)$ is said to be an {\it anisotropy} if it is convex, positively 1-homogeneous (i.e., $F(\lambda x) = \lambda F(x)$ for all $\lambda >0$ and all $x \in \R^N$) and coercive. We will always consider additionally that $F$ is even, that is, a norm. In particular, as all norms in $\R^N$ are equivalent,
there exist  constants $0<c\le C<\infty$ such that
  \begin{equation}\label{FequivEuclid}
  c\Vert \xi\Vert \le F(\xi) \le C \Vert\xi\Vert
  \end{equation}
  where $\Vert \cdot \Vert$ is the Euclidean norm in $\R^N$.

We define the dual or polar function $F^\circ:\R^N\rightarrow [0,\infty)$ of $F$ by
\[F^\circ(\xi):=\sup \{ \xi\cdot\xi^\star \, : \, \xi^\star\in\R^N, \, F(\xi^\star)\le 1\}
=\sup \left\{\frac{\xi \cdot \xi^\star}{F(\xi^\star)}\, : \, \xi^\star\in\R^N\setminus\{0\} \right\}\]
for every $\xi\in\R^N.$
It can be verified that $F^\circ$ is convex, lower semi-continuous and $1$-positively homogeneous. Moreover, \eqref{FequivEuclid} leads to
 \begin{equation}\label{F0equivEuclid}
\frac1{C}\Vert \xi\Vert \le F^\circ(\xi) \le \frac1{c} \Vert\xi\Vert.
\end{equation}

From the definition of $F^\circ$ one gets a Cauchy-Schwarz-type inequality of the form
\begin{equation} \label{CS-ineq}
x \cdot \xi \leq F^\circ(x) F(\xi) \qquad \text{ for all } \ x, \xi \in \R^N.
\end{equation}

The Wulff shape $\mathcal W_R$ of $F$ is defined by
$$\mathcal W_R:=B_{F^\circ}(R):=\{ \xi^\star\in\R^N \, : \, F^\circ(\xi^\star)<R\}.$$
As we are dealing with even anisotropies, $\mathcal W_R$ is a centrally symmetric convex body. We say that $F$ is {\it crystalline} if, furthermore, $\mathcal W_R$ is a convex polytope.




\subsection{Maximal monotone operators on Banach spaces}

Let $X$ be a reflexive Banach space with dual $X'$, and denote by $\big< \,\cdot\,, \,\cdot\,\big>$ the
pairing between $X'$ and $X$.
Let $\mathcal A: X \rightarrow 2^{X'}$ be a multivalued operator on $X$ (equivalently, we write $\mathcal A \subset X \times X'$ for its graph). Hereafter, $D(\mathcal A)$ and $R(\mathcal A)$ mean the domain and range of $\mathcal A$, respectively.

\begin{definition}
$\mathcal A$ is said to be {\it monotone} if
\[\big<\tilde \xi - \tilde\eta, \xi - \eta\big> \geq 0 \qquad \text{for every} \quad  (\xi, \tilde\xi), (\eta, \tilde\eta) \in \mathcal A.\]
Moreover, $\mathcal A$ is called maximal monotone if there exists no other monotone multivalued map whose graph strictly contains the graph of $\mathcal A$.
\end{definition}

Let  $\Phi: X \rightarrow \R \cup \{\infty\}$ be  lower semicontinuous, proper and convex; its subdifferential $\partial \Phi: X \rightarrow 2^{X'}$ is a multivalued operator given as follows:
\begin{equation} \label{orig-def-subD}
\delta \in \partial \Phi(\zeta) \quad \Longleftrightarrow \quad \Phi(\eta) - \Phi(\zeta) \geq \big<\delta, \eta-\zeta\big> \quad \text{for all } \eta \in X.
\end{equation}
This is a well-known example of a maximal monotone operator from $X$ to $X'$ (cf. \cite[Theorem A]{Rock}).

In turn, the duality mapping $\mathcal J_X: X \rightarrow 2^{X'}$ is defined as
\begin{equation} \label{defJ}
\mathcal J_X(u) =\Big\{v \in X' \, : \, \big<u, v\big> = \|u\|^2_X = \|v\|^2_{X'}\Big\}	\qquad \text{for all} \quad u \in X.
\end{equation}

In particular, the
duality mapping $\mathcal J_{L^p}: L^p(\Omega) \rightarrow L^{p'}(\Omega)$ is single-valued and has the following explicit expression:
\begin{equation} \label{defJ-Lp}
\mathcal J_{L^p}(u)(x) = \|u(x)\|^{p-2} u(x) \|u\|_{L^p}^{2-p} \quad \text{for a.e. } x \in \Omega, \quad
\end{equation}
and for all $u \in L^p(\Omega)$.

This will be a crucial characterization (see \cite[Theorem 2.2]{Bar}):
\begin{theorem}[Minty \cite{Minty} and Browder \cite{Brow}]\label{MintyTheorem}
 Let $\mathcal A:  L^p(\Omega) \rightarrow L^{p'}(\Omega)$ be a
monotone functional and let $\mathcal J_{L^{p}}: L^{p}(\Omega) \rightarrow L^{p'}(\Omega)$ be the duality mapping of $L^{p}$. Then $\mathcal A$ is maximal monotone if, and only if, $R(\mathcal A+ \mathcal J_{L^{p}}) = L^{p'}(\Omega)$.
\end{theorem}

Therefore, we need results that ensure that the above range condition holds in order to deduce that an operator is maximal monotone. In particu\-lar, we will apply the following statement (see \cite[Corollary 1.8]{Kinder}):
\begin{theorem}[Hartman and Stampacchia \cite{HS}] \label{Stamp}
Let $\mathfrak K \neq \emptyset$ be a closed convex subset of a reflexive Banach space $X$, and $\mathcal B: \mathfrak K \rightarrow X'$ be monotone, coercive and weakly continuous on finite dimensional subspaces. Then $\mathcal B$ is surjective. 	
\end{theorem}
\noindent Typically the literature refers to the hypotheses in the above theorem as the classical Leray-Lions assumptions (cf. \cite{LL}).

\subsection{Yosida approximation of maps in Hilbert spaces}

Now consider $I$ the identity map on a Hilbert space $H$. Then the {\it resolvent} of $\mathcal A: H \rightarrow 2^H$ is given by
\[J_\lambda =  J_{\mathcal A, \lambda} := (I+\lambda\, \mathcal A)^{-1} \qquad \text{for } \quad \lambda > 0. \]
It holds that the resolvent of a monotone operator $\mathcal A$ is a single-valued non-expansive map from  $R(I + \lambda \, \mathcal A)$ to $H$ (cf. \cite[Proposition 3.5.3]{Au}).

The {\it Yosida approximation} of a maximal monotone map $\mathcal A$ is defined as
$$\mathcal A_\lambda:=\frac{I-J_\lambda}{\lambda},$$
and hence it satisfies
\begin{equation}\label{ApproxSubdifferentialResolvent}
	\mathcal A_\lambda(\xi)\in \mathcal A\left(J_\lambda(\xi)\right) \quad \text{for all} \ \xi \in H.
\end{equation}
In addition, we have the following properties  (see \cite[Theorem 3.5.9]{Au}):
\begin{lemma} \label{prop-Yosida}
{\rm (a)} $\mathcal A_\lambda$ is a single-valued maximal monotone map which is Lipschitz with constant $1/\lambda$.

{\rm (b)} $\mathcal A_\lambda(\xi)$ converges to $\mathcal A^0(\xi)$ as $\lambda \to 0$, where $\mathcal A^0$ denotes the {\it minimal section} of $\mathcal A$ given by
\[\mathcal A^0 (\xi) := \{\eta \in \mathcal A(\xi) \, : \, \eta \text{ has minimal norm in } \mathcal A(\xi)\}.\]
\end{lemma}

Recall that given a proper, convex and lower semicontinuous function $g:\R^N\rightarrow [0, \infty)$, for every $\lambda >0$ the {\it Moreau-Yosida approximation} of $g$ is defined by
\begin{equation}\label{DefMoreauYosida}
	g_\lambda(\xi):=\min_{\eta\in\R^N}\left\{\frac{1}{2\lambda}\Vert \eta - \xi\Vert^2+g(\eta)\right\} \leq g(\xi)
\end{equation}
for all $\xi \in \R^N$. The minimum in \eqref{DefMoreauYosida} is attained at $J_{\partial g,\lambda}(\xi)$, where $\partial g$ denotes the subdifferential of $g$.
Furthermore, $g_\lambda$ is convex and Fr\'echet differentiable with gradient
\begin{equation} \label{grad-MY}
\nabla g_\lambda = (\partial g)_\lambda .
\end{equation}
Thus the Yosida approximation of the subdifferential is equal to the gradient of the Moreau-Yosida approximation.
We also have that $g_\lambda(\xi)\stackrel{\lambda\searrow 0}{\longrightarrow}g(\xi)$ for every $\xi\in\R^N$ (see \cite{Brezis}, \cite[Proposition 1.8]{Show} or \cite[Theorem 6.5.7]{Au}).

In case that $g_\lambda(0) = 0$, the convexity of $g_\lambda$ implies (cf. \cite[Proposition 7.4]{Show})
\begin{equation}\label{CotaInferiorPartialFxixi}
\nabla g_\lambda (\xi)\cdot\xi\ge g_\lambda(\xi) \quad \text{for all} \ \xi \in D(g_\lambda).
\end{equation}

\section{Characterization of the subdifferential on bounded domains} \label{subdif-1}

\subsection{An auxiliary functional containing the subdifferential} 

Hereafter  $\Omega$ will denote an open bounded subset of $\R^N$ with Lipschitz-continuous boundary $\partial \Omega$.
In order to study the minimizers of the functional $\mathcal{F}_{\Omega,\varphi}$ defined by \eqref{defF}, the well-known equivalence \eqref{equiv} tells that one needs to charac\-terize the subdifferential of $\mathcal{F}_{\Omega,\varphi}$ in order to understand the meaning of $0\in\partial\mathcal{F}_{\Omega,\varphi}(u)$.

 Recall that the subdifferential $\partial\mathcal{F}_{\Omega,\varphi}$ of $\mathcal{F}_{\Omega,\varphi}$ is the multivalued operator from $L^p(\Omega)$ to $L^{p'}(\Omega)$ given by
\[v \in \partial\mathcal{F}_{\Omega,\varphi}(u) \quad \Longleftrightarrow\quad \mathcal{F}_{\Omega,\varphi}(w) - \mathcal{F}_{\Omega,\varphi}(u) \geq \int_\Omega v(w-u) \quad \text{for all} \ w \in L^{p}(\Omega).\]

\begin{remark}
It is a well-known result (see e.g. \cite[Proposition 1.6]{Bar}) that the effective domain $\mathcal D(\partial \mathcal F_{\Omega, \varphi}) $ is dense in $\mathcal D(\mathcal F_{\Omega, \varphi})  = W^{1,p}_\varphi(\Omega)$, and consequently in $L^p(\Omega)$. Recall that $u$ belongs to the effective domain of $\mathcal A$ if $\|\mathcal A (u)\|_{X'} < \infty$.
\end{remark}

To get the sought characterization of $\partial\mathcal{F}_{\Omega,\varphi}$, we first introduce an auxiliary functional $\mathcal{A}_\varphi$.

\begin{definition} \label{defA}
	We define $\mathcal{A}_\varphi \subset L^p(\Omega) \times  L^{p'}(\Omega)$ as follows: $(u,v)\in\mathcal{A}_\varphi$ if
	\begin{itemize}
		\item[{\bf (A1)}] $u\in W^{1,p}_\varphi(\Omega)$,
		\item[{\bf (A2)}] $v\in L^{p'}(\Omega)$,
		\item[{\bf (A3)}] there exists $z\in \partial F^p(\nabla u)$ such that $v=-{\rm div} (z)$ weakly, 
meaning that
		$$\int_{\Omega} vw = \int_{\Omega} z\cdot\nabla w \quad \hbox{ for every } w\in W_0^{1,p}(\Omega).$$
		
	\end{itemize}
	
\end{definition}

The goal is to show that both operators $\mathcal{A}_\varphi$ and $\partial \mathcal{F}_{\Omega,\varphi}$ coincide (cf. \autoref{AeqDF}). Once we have established this result, we can conclude that $0 \in \partial \mathcal F_{\Omega,\varphi}(u)$ is equivalent to $u$ satisfying the following:
\begin{equation}\label{EquationBounded}
\left\{ \begin{array}{ll} {\rm div}\big( F^{p-1}(\nabla u) z\big) = 0 & \text{in }  \Omega \smallskip \\
u = \varphi & \text{on }  \partial\Omega \end{array} \right.
\end{equation}
in a weak sense, where $z \in \partial F(\nabla u)$ and we have used the chainrule for subdifferentials  (see \cite[Corollary 16.72]{BC} or \cite[Theorem 2.3.9 (ii)]{Clarke}).

The first inclusion is quite straightforward:
 \begin{lemma}\label{InclusionApartialF}
 $\mathcal{A}_\varphi\subseteq \partial\mathcal{F}_{\Omega,\varphi}.$ In particular, the functional $\mathcal{A}_\varphi$ is monotone.
 \end{lemma}
 \begin{proof}
Let $(u,v)\in \mathcal{A}_\varphi$ and $z$ as in (A3) above. In particular, by definition of subdifferential
\[z \in \partial F^p(\nabla u) \quad \Longrightarrow \quad \int_\Omega  z\cdot\nabla(w-u) \le \int_{\Omega} F^{p}(\nabla w) - \int_{\Omega} F^p(\nabla u) \]
 for every $w\in W^{1,p}_\varphi(\Omega)$. Consequently,
 $$\int_{\Omega} v(w-u) = \int_{\Omega} z\cdot\nabla(w-u) \le \mathcal{F}_{\Omega,\varphi}(w)-\mathcal{F}_{\Omega,\varphi}(u)$$
for every $v \in L^{p'}(\Omega)$. Hence the statement follows by definition of $\partial\mathcal{F}_{\Omega,\varphi}$.
 \end{proof}

\subsection{Properties of the subdifferential of $F^p$} \label{subdif-2}
To establish the remaining inclusion, we need some properties of $\partial F^p$ that will be used repeatedly. We start by writing the elements of $\partial F^p$ in terms of those in $\partial F$ because, as $F$ is  positively 1-homogeneous, its subdifferential has a useful well-known characterization (obtained by choosing $\eta = 0$ and $\eta = 2 \zeta$ in \eqref{orig-def-subD}):
\begin{equation} \label{char-dF}
\delta\in\partial F(\zeta) \quad \Longleftrightarrow \quad \left\{\begin{array}{ll} \delta\cdot\zeta=F(\zeta), &\text{and} \medskip\\ \delta\cdot\eta\le F(\eta) & \text{ for every } \eta\in\R^N. \end{array} \right. 
\end{equation}
Notice that the above inequality is equivalent to $F^\circ(\delta) \leq 1$. In particular, for $\eta = \frac{\delta}{\Vert \delta\Vert}$ it follows by means of \eqref{FequivEuclid} that
\begin{equation} \label{subdif-bdd}
\Vert \delta \Vert \le C \quad  \text{for all} \ \delta\in \partial F(\zeta).
\end{equation}
The latter bound, combined with the chainrule for subdifferentials and again \eqref{FequivEuclid}, leads directly to
\begin{lemma} \label{charFp}
Every $z \in \partial F^p(\xi)$ can be written as
\[z = p F^{p-1}(\xi) \tilde z, \qquad \text{with} \quad \tilde z \in \partial F(\xi)\]
and we can estimate
\[\|z\| \leq p \, C^p \|\xi\|^{p-1}.\]
\end{lemma}

Next we see that the same estimate holds for the corresponding Yosida approximation.
\begin{lemma} \label{aux-Fp}
	The Moreau-Yosida approximation of $F^p$ and the Yosida approximation of $\partial F^p$ satisfy:
	\begin{enumerate}
		\item[{\rm (a)}] If $1 < p<2$, there exists a positive  constant $K_\lambda = K_\lambda(c, p, \lambda)$ with $\displaystyle \lim_{\lambda \searrow 0} K_\lambda = 0$ so that
		\[(F^p)_\lambda(\xi)\ge \frac{c^p}{2} \Vert \xi \Vert^p\quad  \text{if} \quad \Vert \xi\Vert \geq K_\lambda.\]
		\item[{\rm (b)}] $\left\Vert\left(\partial F^p\right)_\lambda(\xi)\right\Vert =\left\Vert\left(\nabla (F^p)_\lambda\right)(\xi)\right\Vert \le p \, C^p \Vert \xi\Vert^{p-1}.$
	\end{enumerate}
	Here $c$ and $C$ are the constants coming from \eqref{FequivEuclid}.
\end{lemma}

 \begin{proof}
 (a) Note that, using \eqref{FequivEuclid} in the definition given by \eqref{DefMoreauYosida}, one has
 $$(F^p)_\lambda(\xi)\ge c^p\min_{\eta\in\R^N}\left\{\frac{1}{2c^p \lambda}\Vert \eta - \xi\Vert^2+\Vert \eta \Vert^p\right\}.$$
 Hence we are done in the case that $\Vert \eta \Vert^p \geq \Vert \xi \Vert^p/2$. Otherwise, if $\Vert \eta \Vert < \Vert \xi \Vert/2^{1/p}$, by the reverse triangle inequality we get
 \[\frac{1}{2c^p \lambda}\Vert \eta - \xi\Vert^2 \geq \frac{1}{2c^p \lambda} (\Vert \xi\Vert - \Vert \eta\Vert)^2 > \frac{\Vert\xi\Vert^2}{2c^p \lambda}(1- 2^{-1/p})^2 > \frac{\Vert \xi \Vert^p}{2},\]
 where the last inequality holds, if and only if,
 \begin{equation} \label{def-Kl}
 \Vert \xi \Vert > \bigg(\frac{\lambda c^p}{(1-2^{-1/p})^2}\bigg)^{\frac1{2-p}} =:K_\lambda.
 \end{equation}

 (b) By \eqref{ApproxSubdifferentialResolvent} with $\mathcal A = \partial F^p$ and the chainrule for subdifferentials, we can write
\begin{equation} \label{chainR}
(\partial F^p)_\lambda(\xi) \in \partial F^p(J_{p,\lambda}(\xi))=pF^{p-1}(J_{p, \lambda }(\xi))\partial F(J_{p, \lambda}(\xi)),
\end{equation}
where $J_{p,\lambda}(\xi) = J_{\partial F^p,\lambda}(\xi)$ satisfies
\[\frac1{2\lambda} \Vert J_{p,\lambda}(\xi)-\xi \Vert^2 + F^p(J_{p,\lambda}(\xi)) \leq \frac1{2\lambda} \Vert \eta-\xi \Vert^2 + F^p(\eta) \quad \text{for all} \ \eta \in \R^N. \]
In particular, for $\eta = \xi$, we get $F^p(J_{p,\lambda}(\xi)) \leq F^p(\xi)$. Raising the inequality to $\frac1{p}$ and using \eqref{FequivEuclid} leads to
\[F(J_{p,\lambda}(\xi)) \leq C \Vert \xi \Vert.\]
This, combined with \eqref{grad-MY} (with $g= F^p$) and \eqref{chainR}, implies
\begin{equation} \label{lemaB-aux}
\left\Vert\left(\partial F^p\right)_\lambda(\xi)\right\Vert =\left\Vert\left(\nabla (F^p)_\lambda\right)(\xi)\right\Vert \le p F^{p-1}(J_{p,\lambda}(\xi)) \Vert \delta\Vert  \leq p C^{p} \Vert \xi\Vert^{p-1},
\end{equation}
for every $\delta \in \partial F(J_{p,\lambda}(\xi))$. Note that we used \eqref{FequivEuclid} and \eqref{subdif-bdd} for the last inequality.
\end{proof}

The existence of minimizers follows from the above properties.
\begin{proposition} \label{Flsc}
There exists at least one function $u$ which is a minimizer of the functional $\mathcal F_{\Omega, \varphi}$ defined in \eqref{defF} for any domain $\Omega \subset \R^N$. 
Moreover, $\mathcal F_{\Omega, \varphi}$ is lower semicontinuous with respect to the weak convergence in $W^{1,p}(\Omega)$.
\end{proposition}

\begin{proof}
As by \eqref{FequivEuclid} the functional $\mathcal F_{\Omega, \varphi}$ is bounded from below, there exists a minimizing sequence $u_n$. By definition of the functional and the Poincaré inequality, $\{u_n\}_{n\in \N}$ is uniformly bounded in $W^{1,p}(\Omega)$ and hence it converges subsequentially and weakly in $W^{1,p}(\Omega)$ to $u \in W^{1,p}(\Omega)$.

Now let $z \in \partial F^p(\nabla u)$, then $z \in L^{p'}(\Omega)$ and $z = p F^{p-1}(\nabla u) \delta$ with $\delta \in \partial F(\nabla u)$ by Lemma \ref{charFp}. By \eqref{char-dF} and the weak convergence, we have
\begin{align*}
 \int_\Omega F^p(\nabla u) & = \frac1{p}\int_\Omega z \cdot \nabla u = \frac1{p} \lim_{n \to \infty} \int_\Omega z \cdot \nabla u_n \leq \lim_{n \to \infty}  \int_\Omega F^{p-1}(\nabla u) F(\nabla u_n) \\
 & \leq  \Big(\int_\Omega F^p(\nabla u)\Big)^{\frac1{p'}}  \liminf_{n\to \infty} \Big(\int_\Omega F^p(\nabla u_n)\Big)^{\frac1{p}},
\end{align*}
where the latter follows by application of Hölder's inequality. Thus
\[ \int_\Omega F^p(\nabla u) \leq \liminf_{n\to \infty} \int_\Omega F^p(\nabla u_n),\]
as desired.
\end{proof}
\noindent It remains to show that such a minimizer is a weak solution of \eqref{EquationBounded}, which will be a direct consequence of the characterization of $\partial \mathcal F_{\Omega, \varphi}$.

Throughout this paper, we will perform several approximation arguments whose  outcome is a suitable weak limit $z \in L^{p'}(\Omega)$. The following result, in the spirit of the {\it Minty-Browder technique} (see \cite{Min}),  gives a practical criterion to guarantee that $z$ is indeed in the appropriate subdifferential. For brevity, here and in the sequel we will use the notation
\begin{equation} \label{defPsi}
\Psi(\zeta) = \|\zeta\|^{p-2} \zeta.
\end{equation}
\begin{lemma} \label{Minty-tec}
Given $u \in W^{1,p}(\Omega)$, $z \in L^{p'}(\Omega)$ and $\alpha \in\R$, suppose that the following inequality holds for a.e. $x\in \Omega$ and every $\xi\in\R^N$:
	\begin{equation} \label{Minty2}
	\big(z(x) -z'
	-\alpha \Psi(\xi)\big)\cdot(\nabla u (x) - \xi) \ge 0 \quad \hbox{for every } z'\in\partial F^p(\xi).
	\end{equation}
Then one can conclude that $z \in \partial F^p(\nabla u)+ \alpha \Psi(\nabla u)$. The same conclusion holds if we only require that \eqref{Minty2} holds for $z' = (\partial F^p)^0(\xi)$.
\end{lemma}

\begin{proof}
Notice that it is enough to prove the second claim, but we keep the stronger hypothesis \eqref{Minty2} because of the applications afterwards.	Therefore, there exists a null set $N_1\subset\R^N$ such that
	$$\big(z(x)-(\partial F^p)^0(\xi)-\alpha
	\Psi(\xi)\big)\cdot (\nabla u(x)-\xi)\ge 0$$
	for every $x\in \Omega\setminus N_1$ and $\xi\in\R^N\setminus N_x$,
	where $N_x\subset\R^N$ is null. Let $x\in\Omega\setminus N_1$. Then, for a.e. $0\neq\eta\in \R^N$, $\{t \, : \, \nabla u(x) + t\eta \in N_x\}$ is null and we may take a sequence $\{\varepsilon_n\}_{n\in\N}$ with $\varepsilon_n \searrow 0 $ and $\{\nabla u(x) + \varepsilon_n \eta\}_{n\in\N}\subset \R^N\setminus N_x$. If we are in this case, 
	 we can write
	\begin{align*}
	 \big(z(x)&-(\partial F^p)^0(\nabla u(x) + \varepsilon_n \eta) - \alpha \Psi(\nabla u(x) + \varepsilon_n \eta)\big)\cdot\eta \\ & \quad =\frac{1}{\varepsilon_n}\left(z(x)-(\partial F^p)^0(\nabla u(x) + \varepsilon_n \eta) - \alpha \Psi(\nabla u(x) + \varepsilon_n \eta)\right)\cdot(\varepsilon_n \eta)  \le 0 \quad
	 \end{align*}
for every $n\in\N$.
	
	Now, since $(\partial F^p)^0(\nabla u(x) + \varepsilon_n \eta)$ is bounded we may suppose that, up to a subsequence, $(\partial F^p)^0(\nabla u(x) + \varepsilon_n \eta)\rightarrow \tilde{z}\in \partial F^p(\nabla u(x))$ as $n\to\infty$ (recall that $\partial F^p$ is closed since it is maximal monotone). Hence, by the continuity of $\Psi$, taking limits in the previous equation we get that
	$$\Big(z(x)-\tilde{z} -\alpha \Psi(\nabla u(x))\Big)\cdot \eta \le 0.$$
	Since this holds for a.e. $\eta\in\R^N$, we conclude that $z(x)-\alpha \Psi(\nabla u(x))=\tilde{z}\in\partial F^p(\nabla u(x))$.
\end{proof}

\subsection{A range condition for the auxiliary functional}

 To prove Theo\-rem \ref{AeqDF} we need to show that the operator $\mathcal{A}_\varphi$ introduced in Definition \ref{defA} satisfies a range condition, which involves the duality mapping of $L^p(\Omega)$ given by \eqref{defJ-Lp}.
 \begin{proposition}\label{Range-0}
If $\varphi\in W^{1-\frac1{p}, p}(\partial\Omega)$, then $L^{p'}(\Omega)= R(\mathcal J_{L^p}+\mathcal{A}_\varphi)$.

 \end{proposition}

 \begin{proof}
 Given $f\in L^{p'}(\Omega)$, we aim to prove that there exists $u\in W^{1,p}_\varphi(\Omega)$ such that $(u, f-\mathcal J_{L^p} u)\in\mathcal{A}_\varphi$. With this purpose, let us consider  \[\mathcal{A}_{\varphi}^{n,m}:  W^{1,p}_\varphi(\Omega) \longrightarrow (W^{1,p}(\Omega))',\]
  which is a sequence of operators approximating $\mathcal{A}_\varphi$ and defined as follows:
 \begin{equation} \label{def-aproxA}
 v= \mathcal A_\varphi^{n,m}(u)= -\displaystyle \hbox{div}\left(\big(\nabla (F^p)_\frac{1}{n}\big)(\nabla u) + \frac1{m} \|\nabla u\|^{p-2} \nabla u\right)
 \end{equation}
  in the weak  sense. In particular, with the notation from \eqref{defPsi} it holds
 	$$\int_{\Omega}vw=\int_{\Omega}\big(\nabla(F^p)_\frac{1}{n}\big)(\nabla u)\cdot\nabla w + \frac1m \int_{\Omega} \Psi(\nabla u) \cdot \nabla w
 	\quad \text{for every } \ w\in W_0^{1,p}(\Omega).	$$

{\bf Claim 1.} {\it The operators $\mathcal{A}_{\varphi}^{n,m}+ \mathcal J_{L^p}$ satisfy the hypotheses of \autoref{Stamp}}.
\begin{proof}[Proof of Claim 1]
First, $\mathcal J_{L^p}$ is uniformly continuous on every bounded subset of $L^p(\Omega)$ (cf. \cite[Theorem 1.2]{Bar}) and it is monotone because it coincides with the subdifferential of $\frac1{2} \|\,\cdot\|^2_{L^p}$ (cf. \cite[p.7]{Bar}).

Since $\mathcal{A}_{\varphi}^{n,m}(u) = -\displaystyle \hbox{div}\left(\left(\nabla (F^p)_\frac{1}{n}\right)(\nabla u)\right) -\frac1{m} \Delta_p u$,  monotonicity and continuity for the first term follow from \eqref{grad-MY} and Lemma \ref{prop-Yosida} (a). In turn, the $p$-Laplacian is known (cf. \cite[Lemma 2.2]{JS}) to be Hölder continuous on  $W^{1,p}(\Omega)$, and the monotonicity follows from the inequalities (see \cite[(2.10)]{JS})
\[\big(\Psi(\xi) - \Psi(\eta)\big) \cdot (\xi-\eta) \geq \left\{\begin{array}{ll} c_p \|\xi-\eta\|^p, & \text{if } p \geq 2, \medskip \\
\displaystyle c_p \frac{\|\xi -\eta\|^2}{\big(\|\xi\|+\|\eta\|\big)^{2-p}} & \text{if } 1<p<2,\end{array}
\right.\]
which hold for every $\xi, \eta \in \R^N$. This, combined with Poincaré inequality, easily leads to coercivity when $p\geq 2$; in fact, for any $u, w \in W^{1,p}_\varphi(\Omega)$ we get
\[\frac{\big<\Psi(\nabla u)-\Psi(\nabla w), \nabla u-\nabla w\big>}{\|u - w\|_{W^{1,p}} } \geq C_{p, \Omega} \|\nabla u -\nabla w\|^{p-1}_{L^p} \to \infty \ \text{as} \ \|\nabla u\|_{L^p} \to \infty.\]
To establish the coercivity of  $\mathcal{A}_{\varphi}^{n,m}+ \mathcal J_{L^p}$ in the remaining case $1 < p<2$, again by Poincaré inequality, it is enough to show that, given a fixed $w \in W^{1,p}_\varphi(\Omega)$, it holds
\[\frac{\big<\Psi(\nabla u)-\Psi(\nabla w), \nabla u-\nabla w\big>}{\|\nabla u -\nabla w\|_{L^p} } \to \infty \quad \text{as} \quad \|\nabla u\|_{L^p} \to \infty\]
for any $u \in W^{1,p}_\varphi(\Omega)$. With this aim, the triangle inequality leads to
\begin{align*}
\frac{\big<\Psi(\nabla u), \nabla u\big> + \big<\Psi(\nabla w), \nabla w\big>}{\|\nabla u -\nabla w\|_{L^p} } \geq \frac{\|\nabla u\|_{L^p}^p + \|\nabla w\|_{L^p}^p}{\|\nabla u\|_{L^p} + \|\nabla w\|_{L^p}},
\end{align*}
which tends to $\infty$ as $\|\nabla u\|_{L^p} \to \infty$.
It remains to check that the crossed terms keep bounded in the limit. Indeed, we can assume that $\|\nabla u\|_{L^p} > \|\nabla w\|_{L^p}$, and then
\begin{align*}
\bigg|\frac{\big<\Psi(\nabla u), \nabla w\big> + \big<\Psi(\nabla w), \nabla u\big>}{\|\nabla u -\nabla w\|_{L^p}}\bigg| & \leq \frac{\displaystyle\int_\Omega \|\nabla u\|^{p-1} \|\nabla w\|}{\|\nabla u -\nabla w\|_{L^p}} + \frac{\displaystyle \int_\Omega \|\nabla w\|^{p-1} \|\nabla u\|}{\|\nabla u -\nabla w\|_{L^p}} \\ & \leq \frac{\|\nabla u\|_{L^p}^{p-1} \|\nabla w\|_{L^p}}{\|\nabla u\|_{L^p}-\|\nabla w\|_{L^p}} + \frac{\|\nabla w\|_{L^p}^{p-1} \|\nabla u\|_{L^p}}{\|\nabla u\|_{L^p}-\|\nabla w\|_{L^p}},
\end{align*}
where we applied Hölder and the reverse triangle inequality. Notice that, as $p \in (1,2)$, the right hand side tends to $0 + \|\nabla w\|_{L^p}^{p-1}$ as  $\|\nabla u\|_{L^p} \to \infty$, which is a bounded limit, as desired.
\end{proof}

 Now, since $ W^{1,p}_\varphi(\Omega)$ is closed and convex in $W^{1,p}(\Omega)$,  by \autoref{Stamp} we deduce that
\begin{equation} \label{Leray}
L^{p'}(\Omega)\subset (W^{1,p}(\Omega))'\subseteq R(\mathcal J_{L^p}+\mathcal{A}_{\varphi}^{n,m}).
\end{equation}
Thus there exist $u_{n,m}\in W^{1,p}_\varphi(\Omega)$ such that
\begin{equation} \label{AfJ}
 \mathcal{A}_{\varphi}^{n,m}(u_{n,m}) =f-\mathcal J_{L^p}(u_{n,m}) \in L^{p'}(\Omega).
 \end{equation}
Therefore,
\begin{equation}\label{equationAn}
\int_{\Omega}\!\!(f- \mathcal J_{L^p}(u_{n,m}))w=\int_{\Omega}\!\!\big(\nabla(F^p)_\frac{1}{n}\big)(\nabla u_{n,m})\cdot\nabla w + \frac1m \int_{\Omega}\!\! \Psi(\nabla u_{n, m}) \cdot \nabla w
\end{equation}
for every  $w\in W_0^{1,p}(\Omega)$.

{\bf Claim 2.} {\it $\{u_{n,m}\}_{n\in\N}$ is uniformly bounded in $W^{1,p}(\Omega)$.}

\begin{proof}[Proof of Claim 2]
Given  $w\in W^{1,p}_\varphi(\Omega)$, letting $w-u_{n,m}$ as a test function in \eqref{equationAn} yields
\begin{align*}
\int_{\Omega}\big(f-\mathcal J_{L^p}(u_{n,m})\big)(w-u_{n,m}) = &\int_{\Omega}\big(\nabla(F^p)_\frac{1}{n}\big)(\nabla u_{n,m})\cdot (\nabla w- \nabla u_{n,m}) \\ & + \frac1m \int_{\Omega} \Psi(\nabla u_{n,m}) \cdot \nabla w - \frac1{m} \|\nabla u_{n,m}\|_{L^p}^p.
\end{align*}
First, the integral on the left hand side, by means of the definition of the duality mapping in \eqref{defJ}, can be written as
\begin{align*}
	\int_{\Omega}\!\!\big(f-\mathcal J_{L^p}(u_{n,m})\big)(w-u_{n,m}) = \int_\Omega \!\! f (w - u_{n,m}) \! -\! \int_\Omega \!\! \mathcal J_{L^p}(u_{n,m}) w  + \|u_{n,m}\|_{L^p}^2.
\end{align*}
Rearranging terms, we obtain
\begin{align*}
 \|u_{n,m}\|_{L^p}^2 & + \frac1{m} \|\nabla u_{n,m}\|_{L^p}^p + \int_{\Omega} \!\!\big(\nabla(F^p)_\frac{1}{n}\big)(\nabla u_{n,m})\cdot\nabla u_{n,m} = \int_\Omega \!\! f (u_{n,m} -w) +  \\ & + \int_\Omega \!\!\mathcal J_{L^p}(u_{n,m}) w  +
\int_{\Omega}  \!\!\big(\nabla(F^p)_\frac{1}{n}\big)(\nabla u_{n,m})\!\cdot\!\nabla w +  \frac1m \int_{\Omega}  \!\!\Psi(\nabla u_{n,m})\! \cdot \! \nabla w.
\end{align*}
Now, by \eqref{CotaInferiorPartialFxixi} for $g=F^p$ and by H\"older's inequality, we can estimate
\begin{align*}
	\|u_{n,m}\|_{L^p}^2 + \frac{\|\nabla u_{n,m}\|_{L^p}^p}{m} & \leq  \|f\|_{L^{p'}} \big(\|w\|_{L^p} + \|u_{n,m}\|_{L^p} \big) \!+ \!\|\mathcal J_{L^p}(u_{n,m})\|_{L^{p'}} \|w\|_{L^p}  \\ & + \!\!
	\int_{\Omega}\big(\nabla(F^p)_\frac{1}{n}\big)(\nabla u_{n,m})\cdot\nabla w +  \frac{\|\nabla u_{n,m}\|^{p-1}_{L^p}}{m}  \|\nabla w\|_{L^p}.
\end{align*}

On the other hand, by Cauchy-Schwarz for the inner product in $\R^N$, Lemma \ref{aux-Fp} (b) and Young's inequality $\Big[ab\leq \varepsilon a^{p'}+ C(\varepsilon) b^p$, where $C(\varepsilon) = \frac{(\varepsilon p')^{1-p}}{p}\Big]$ with $\varepsilon = \frac1{2m p C^p}$, 
we have that
\begin{align*}\int_{\Omega}\big(\nabla(F^p)_\frac{1}{n}\big)(\nabla u_{n,m})\cdot\nabla w
& \le p \, C^{p} \int_{\Omega}\Vert \nabla u_{n,m}\Vert^{p-1}\Vert\nabla w\Vert
\\ & \le\frac{1}{2 m}\int_{\Omega}\Vert \nabla u_{n,m}\Vert^p+ 
m^{p-1} \mathcal C\int_{\Omega}\Vert \nabla w\Vert^p,
\end{align*}
where $\mathcal C$ denotes hereafter any constant which depends only on $p$ and the constants from \eqref{FequivEuclid}, whose concrete meaning may change from line to line (and we use $\tilde{\mathcal C}$ whenever two of such constants appear on the same line).

Bringing these bounds together, using $\|\mathcal J_{L^p}(u_{n,m})\|_{L^{p'}} = \|u_{n,m}\|_{L^{p}}$, Cauchy's inequality $a b \leq \frac{a^2}{2}+\frac{b^2}{2}$ and  Young's inequality with $\varepsilon = 1/4$, we get
\begin{align*}
	\|u_{n,m}\|_{L^p}^2 + \frac{1}{2    m} \|\nabla u_{n,m}\|_{L^p}^p  &\leq  \|f\|_{L^{p'}} \|w\|_{L^p}  +   \frac1{2}\|u_{n,m}\|_{L^p}^2 + \frac1{2} \big(\|f\|_{L^{p'}}+ \|w\|_{L^p}\big)^2  \nonumber \\ & \hspace*{-0.7cm} +
	m^{p-1} \mathcal C \Vert \nabla w\Vert^p_{L^p} +  \frac1{4 m}  \|\nabla u_{n,m}\|^{p}_{L^p} + \left(\frac{4}{p'}\right)^{p-1}\frac{ \|\nabla w\|_{L^p}^p}{pm}.
\end{align*}
In short, we reach
\[
\frac1{2} \|u_{n,m}\|_{L^p}^2 + \frac1{4 m} \|\nabla u_{n,m}\|_{L^p}^p  \leq  \frac3{2} \big(\|f\|_{L^{p'}}^2 + \|w\|_{L^p}^2\big) + (m^{p-1} + 1) \mathcal C \Vert \nabla w\Vert^p_{L^p},	
\]
from which the claim follows.
\end{proof}

Accordingly, up to a subsequence, we may assume that $u_{n,m}$ converges in $L^p(\Omega)$ and a.e. to some $u_m\in L^p(\Omega)$ as $n\to\infty$; and $\nabla u_{n,m} \rightharpoonup \nabla u_m$ weakly in $L^{p}(\Omega)$ as $n\to\infty$.  The continuity of the trace operator with respect to the weak convergence in $W^{1,p}(\Omega)$ (see \cite[Corollary 18.4]{Leoni}) guarantees that $u_m \in W^{1,p}_\varphi(\Omega)$.

{\bf Claim 3.} {\it There exists $z_m \in \partial F^p(\nabla u_m)$ such that
\begin{equation} \label{test-eqn3}
	\int_{\Omega}(f- \mathcal J_{L^p}(u_m))\psi=\int_{\Omega}z_m\cdot\nabla \psi + \frac1{m} \int_\Omega \Psi(\nabla u_m) \cdot \nabla \psi,
\end{equation}
for every $\psi\in W^{1,p}_0(\Omega)$.}

\begin{proof}[Proof of Claim 3]
First, Lemma \ref{aux-Fp} (b) ensures that
$$\left\Vert\left(\nabla (F^p)_\frac{1}{n}\right)(\nabla u_{n,m})\right\Vert \le pC^{p} \Vert\nabla u_{n,m}\Vert^{p-1}$$
thus
$$\left\{\left(\nabla (F^p)_\frac{1}{n}\right)(\nabla u_{n,m})\right\}_{n\in\N} \quad \text{is bounded in } \  L^{p'}(\Omega).$$
 Consequently, we may assume that
\begin{equation}\label{ConvToz}
\left(\nabla (F^p)_\frac{1}{n}\right)(\nabla u_{n,m})\rightharpoonup \overline z_m  \ \hbox{ weakly in  } L^{p'}(\Omega) \ \hbox{ as } n\to\infty.\end{equation}
In turn, we also have that $\|\Psi(\nabla u_{n,m})\| =\Vert\nabla u_{n,m}\Vert^{p-1}$, and hence
\begin{equation}\label{ConvTog}
	\Psi(\nabla u_{n,m})\rightharpoonup g_m \ \hbox{ weakly in  } L^{p'}(\Omega) \ \hbox{ as } n\to\infty .\end{equation}
Therefore, taking $\psi\in W^{1,p}_0(\Omega)$ as test function in \eqref{equationAn} and letting $n\to\infty$ we get 
\begin{equation} \label{test-eqn2}
\int_{\Omega}(f- \mathcal J_{L^p}(u_m))\psi=\int_{\Omega}\overline{z}_m\cdot\nabla \psi + \frac1{m} \int_\Omega g_m \cdot \nabla \psi,
\end{equation}
in particular,
$$-{\rm div} \Big(\overline{z}_m + \frac1{m} g_m\Big)=f- \mathcal J_{L^p}(u_m) \ \ \hbox{in the weak sense}$$
and
\begin{equation}\label{ConvToDivz}
	\mathcal{A}_{\varphi}^{n,m}(u_{n,m}) \rightharpoonup - {\rm div}\left(\overline{z}_m +\frac{1}{m}g_m \right) \  \text{ weakly in $L^{p'}(\Omega)$} \ \text{as $n\to\infty$}, \  \end{equation}
 which comes from \eqref{AfJ}.  

Moreover, since for a fixed $\xi\in\R^N$,  $\big(\nabla(F^p)_\frac{1}{n}\big)(\xi)=\left(\partial F^p\right)_{\frac{1}{n}}(\xi)$ tends to $(\partial F^p)^0(\xi)$ as $n\to\infty$ (by Lemma \ref{prop-Yosida} (b)) and $\left\Vert\big(\nabla(F^p)_\frac{1}{n}\big)(\xi)\right\Vert\le pCF^{p-1}(\xi)$ (by Lemma \ref{aux-Fp} (b)), by the dominated convergence theorem, we have that
\begin{equation} \label{ConvOfNablaFp1n}
\left(\nabla(F^p)_\frac{1}{n}\right)(\nabla g) \rightarrow  (\partial F^p)^0(\nabla g)\ \ \hbox{ in } L^{p'}(\Omega)\hbox{ as } n\to\infty
\end{equation}
for every $g\in W^{1,p}(\Omega)$. 

Next we aim to show that $\overline z_m + \frac1{m} g_m\in\partial F^p(\nabla u_m)+ \frac1{m} \Psi(\nabla u_m)$ by means of Lemma \ref{Minty-tec}. 
With this aim, let $\xi\in\R^N$ and take $\omega:\Omega\rightarrow \R$ defined by $\omega(\eta):=\xi\cdot\eta$ so that $\nabla \omega =\xi$. Then, letting $\phi\in C^\infty_0(\Omega)$ such that $\phi\ge 0$, by the monotonicity of $\nabla(F^p)_\frac{1}{n}=\left(\partial F^p\right)_{\frac{1}{n}}$ and $\Psi$ we have, for every $n\in\N$,
$$\begin{array}{l}\displaystyle 0\le \!\big<\big(\nabla(F^p)_\frac{1}{n}\big)(\nabla u_{n,m})-\big(\nabla(F^p)_\frac{1}{n}\big)(\nabla \omega), \nabla(u_{n,m}-\omega)\phi\big>\\ [12pt]
\displaystyle \phantom{0\le} + \left<\frac{\Psi(\nabla u_{n,m}) -\Psi(\nabla \omega)}{m}, \nabla(u_{n,m}-\omega)\phi\right> \\ [12pt]
\phantom{0}= \displaystyle \!\int_{\Omega} 
\!\!\mathcal{A}_{\varphi}^{n,m}(u_{n,m})(u_{n,m}-\omega)\phi\\ [12pt]
\displaystyle \phantom{0\le}-\int_{\Omega}\left(\big(\nabla(F^p)_\frac{1}{n}\big)(\nabla u_{n,m})+\frac{1}{m}\Psi(\nabla u_{n,m})\right)\cdot\nabla\phi(u_{n,m}-\omega)\\ [12pt]
\displaystyle \phantom{0=}-\int_{\Omega}\left(\big(\nabla(F^p)_\frac{1}{n}\big)(\nabla \omega)+\frac{1}{m}\Psi(\nabla \omega)\right)\cdot (\nabla u_{n,m}-\nabla \omega)\phi
\end{array}$$
thus,  letting $n\to\infty$, by \eqref{ConvToz}, \eqref{ConvTog}, \eqref{ConvToDivz} and \eqref{ConvOfNablaFp1n} we get
\begin{align*}
0 \leq & - \int_\Omega {\rm div}\Big(\overline{z}_m + \frac1{m} g_m\Big)(u_m - \omega) \phi  - \int_\Omega \Big(\overline{z}_m + \frac1{m}g_m\Big) \cdot \nabla \phi (u_m - \omega) \\ & - \int_{\Omega} \Big((\partial F^p)^0(\nabla \omega) + \frac1{m} \Psi(\nabla \omega)\Big)\cdot(\nabla u_m-\nabla \omega) \phi \\ & \!\!= \int_\Omega \Big(\overline{z}_m + \frac1{m} g_m - (\partial F^p)^0(\nabla \omega)- \frac1{m} \Psi(\nabla \omega)\Big) \cdot(\nabla u_m-\nabla \omega) \phi.
\end{align*}
Hence, since $\phi\in C^\infty_0(\Omega)$ is arbitrary and $\nabla \omega = \xi$, we conclude \eqref{Minty2} for  $z' = (\partial F^p)^0(\xi)$ and $\alpha = \frac1{m}$. Then $\overline{z}_m+\frac{1}{m}g_m \in\partial F^p(\nabla u_m)+ \frac1{m}\Psi(\nabla u_m)$  by application of Lemma \ref{Minty-tec}. By substitution in \eqref{test-eqn2}, the claim follows.
\end{proof}

{\bf Claim 4.} {\it $\{u_m\}_{m\in \mathbb N}$ is uniformly bounded in $W^{1,p}(\Omega)$.}
\begin{proof}[Proof of Claim 4]
  We proceed as before, but now taking $w - u_m$ as test function in \eqref{test-eqn3} with $w\in W^{1,p}_\varphi(\Omega)$. Then, as $z_m = p F^{p-1}(\nabla u_m) \tilde z_m$ with $\tilde z_m \in \partial F(\nabla u_m)$,  by means of \eqref{char-dF} we can estimate
\begin{align*}
	\int_{\Omega}(f- \mathcal J_{L^p}(u_m))(w - u_m) \leq & \, p \int_{\Omega} F^{p-1}(\nabla u_m) F(\nabla w) - p  \int_{\Omega} F^p(\nabla u_m) \\ &  + \frac1{m} \int_\Omega \Psi(\nabla u_m) \cdot  \nabla w  - \frac1{m} \|\nabla u_m\|_{L^p}^p.
\end{align*}
Rearranging terms, using \eqref{defJ} and applying Hölder's inequality, we have
\begin{align*}
	\|u_{m}\|_{L^p}^2 + & \left(p \, c^p +\frac1{m}\right)  \|\nabla u_{m}\|_{L^p}^p  \leq  \|f\|_{L^{p'}} \big(\|w\|_{L^p} + \|u_{m}\|_{L^p} \big) \\ &  \hspace*{0.7cm} + \|\mathcal J_{L^p}(u_{m})\|_{L^{p'}} \|w\|_{L^p}    +\left(\frac1m + p C^p\right)  \|\nabla u_{m}\|^{p-1}_{L^p}  \|\nabla w\|_{L^p},
\end{align*}
where $c, C$ are the constants coming from \eqref{char-dF}. By definition of the duality map, Cauchy's inequality and Young's inequality with $\varepsilon = \frac1{2} \frac{p c^p + 1/m}{p C^p + 1/m} $, 
\begin{align*}
	\frac1{2}\|u_{m}\|_{L^p}^2 + \frac1{2}\Big(p \, c^p +\frac1{m}\Big) \|\nabla u_{m}\|_{L^p}^p & \leq  \frac{3}{2} \big(\|f\|_{L^{p'}}^2 + \| w\|_{L^p}^2\big)+ C(\varepsilon)\|\nabla w\|_{L^p}^p,
\end{align*}
where $C(\varepsilon) = \frac1{p}\Big(\frac{2 (p C^p +1/m)}{ p'(p c^p +1/m)}\Big)^{p-1} \leq \frac1{p}\Big(\frac{2(p C^p +1)}{p'p c^p}\Big)^{p-1} $. Accordingly,  we reach
\begin{align*}
	\frac1{2}\|u_{m}\|_{L^p}^2 + \mathcal C \|\nabla u_{m}\|_{L^p}^p & \leq  \tilde{\mathcal C} \big(\|f\|_{L^{p'}}^2 + \| w\|_{L^p}^2+ \|\nabla w\|_{L^p}^p\big),
\end{align*}
as desired.\end{proof}

Therefore, $u_m$ converges subsequentially in $L^p(\Omega)$ and a.e. to some $u  \in L^p(\Omega)$ as $m \to \infty$; and $\nabla u_m \rightharpoonup \nabla u$ weakly in $L^p(\Omega)$ as $m \to \infty$.  The continuity of the trace operator with respect to the weak convergence in $W^{1,p}(\Omega)$ ensures that $u \in W^{1,p}_\varphi(\Omega)$.

 Next, by Lemma \ref{charFp}, as $z_m \in \partial F^p(\nabla u_m)$, we have
\[\|z_m\| \leq p \,  C^{p-1} \|\nabla u_m\|^{p-1}.\]
This implies that $\{z_m\}$ is uniformly bounded in $L^{p'}(\Omega)$, and hence $z_m \rightharpoonup z$ weakly in $L^{p'}(\Omega)$ as $m \to \infty$. On the other hand, by Claim 4, $\{\Psi(\nabla u_m)\}_{m\in \mathbb N}$ is bounded in $L^{p'}(\Omega)$ so
\begin{equation} \label{Psito0}
\frac{1}{m}\Psi(\nabla u_m)\rightharpoonup 0 \  \hbox{ weakly in  } L^{p'}(\Omega) \ \hbox{ as } m\to\infty.
\end{equation}
Thus, taking limits as $m\to\infty$ in \eqref{test-eqn3}, we get that
\[\int_{\Omega}(f- \mathcal J_{L^p}(u))\psi=\int_{\Omega}z\cdot\nabla \psi,
 \ \ \hbox{ for every } \psi\in W_0^{1,p}(\Omega),\]
that is,
\begin{equation} \label{divz}
	-{\rm div} (z)=f- \mathcal J_{L^p}(u) \ \ \hbox{in the weak sense}.
\end{equation}

It remains to see that $z\in\partial F^p(\nabla u)$, which  can be obtained after establishing \eqref{Minty2}. 
Therefore, let $\xi\in\R^N$, $z'\in \partial F^p(\xi)$ and take $\omega:\Omega\rightarrow \R$ defined by $\omega(\eta):=\xi\cdot\eta$ so that $\nabla \omega =\xi$. Then, letting $\phi\in C^\infty_0(\Omega)$ such that $\phi\ge 0$, by the monotonicity of $\partial F^p + \frac1{m}\Psi$ and Claim 3 we have, for every $m\in\N$,
$$\begin{array}{l}\displaystyle 0\le\int_{\Omega}\left(z_m+\frac1{m} \Psi(\nabla u_m) -z' -\frac1{m}\Psi(\nabla \omega)\right)\cdot (\nabla u_m-\nabla \omega)\phi   \\ [12pt]
	\,\,\,\,= \displaystyle \phantom{0} \int_{\Omega}(f- \mathcal J_{L^p}(u_m)) (u_m-\omega)\phi -\int_{\Omega}\Big(z_m+\frac1{m} \Psi(\nabla  u_m)\Big)\cdot\nabla\phi(u_m-\omega)  \\ [12pt] \displaystyle  	\,\,\,\ -\int_{\Omega}\left(z' + \frac1{m} \Psi(\nabla \omega)\right)\cdot (\nabla u_m-\nabla \omega)\phi
\end{array}$$
thus, since by Fatou's Lemma
$$\liminf\limits_{n\to\infty} \int_{\Omega} \mathcal J_{L^p}(u_m) \cdot u_m\ge \int_{\Omega}  \mathcal J_{L^p}(u) \cdot u,$$
taking the infimum limit as $m\to\infty$, by means of \eqref{Psito0} we get
$$\begin{array}{l}\displaystyle 0\le \int_{\Omega}\big(f-\mathcal J_{L^p}(u)\big) (u-\omega)\phi-\int_{\Omega}z\cdot\nabla\phi(u-\omega)-\int_{\Omega}z'\cdot (\nabla u-\nabla \omega)\phi \\ [14pt]
	\displaystyle \phantom{0}=\int_{\Omega} \left(z-z'\right)\cdot(\nabla u-\nabla\omega)\phi,
\end{array}$$
where the last equality comes from \eqref{divz}.
Consequently, since $\phi\in C^\infty_0(\Omega)$ is arbitrary and $\nabla \omega = \xi$, we conclude \eqref{Minty2} with $\alpha =0$, as desired. 
Finally, Lemma \ref{Minty-tec} yields $z\in\partial F^p(\nabla u)$, which leads to $(u,f- \mathcal J_{L^p}(u))\in\mathcal{A}_\varphi$, which finishes the proof.
\end{proof}

\begin{remark}
Let us point out that in the singular case, that is, when $1<p<2$, there is a significative shortcut of the above proof, as one can exploit Lemma \ref{aux-Fp} (a) to prove that the operators $\mathcal A_\varphi^n(u):= -\displaystyle \hbox{\rm div}\left(\big(\nabla (F^p)_\frac{1}{n}\big)(\nabla u)\right)$ are coercive, and then the result follows by performing only one limiting argument. The drawback of this shorter proof is that some bounds depend explicitly on the measure of $\Omega$, and in the next section we plan to extend the results (with minor changes) to unbounded domains.
\end{remark}

\subsection{Proof of the characterization} After all the machinery developed in the previous subsections, the following result is almost straightforward.
\begin{theorem} \label{AeqDF}
	$\mathcal{A}_\varphi=\partial \mathcal{F}_{\Omega,\varphi}$ for any $\varphi \in W^{1-\frac1{p},p}(\partial \Omega)$.
\end{theorem}
\begin{proof}
Recall that, by Lemma \ref{InclusionApartialF}, $\mathcal{A}_\varphi\subseteq \partial\mathcal{F}_{\Omega,\varphi}$. Next by Proposition \ref{Range-0} we have that
$$L^{p'}(\Omega)=R\left(\mathcal J_{L^p}+\mathcal{A}_\varphi\right) $$
thus, by Theorem \ref{MintyTheorem}, we have that $\mathcal{A}_\varphi$ is a maximal monotone operator, and hence the claim follows.\end{proof}
\noindent By Definition \ref{defA}, Lemma \ref{charFp} and \eqref{char-dF}, this finishes the proof of \autoref{Emin}.

\section{Characterization of the subdifferential in unbounded exterior domains} \label{unbdd-sec}


After suitable adaptations that we will sketch here, we can reproduce the proof of the previous section and make it work for unbounded domains of the form $\mathbb R^N \setminus \overline \Omega$, where $\Omega$ is a bounded domain with Lipschitz boundary. The main difference is that we work with the homogeneous Sobolev space $\dot W^{1,p}(\R^N)$, defined as the completion of $C^\infty_c(\R^N)$ with respect to the norm
\begin{equation} \label{norm-dot}
\|u\|_{\dot W^{1,p}}:= \|\nabla u\|_{L^p},
\end{equation}
see e.g.~\cite{Bras}. Setting $p^\ast = \frac{Np}{N-p}$, $\dot W^{1,p}(\R^N)$ can be identified with the following space of functions:
\[\left\{u \in L^{p^\ast}(\R^N) \, : \, \nabla u \in L^p\left(\R^N;\R^N\right)\right\}.\] By means of the Gagliardo–Nirenberg–Sobolev inequality, this ensures that
\[\|u\|_{L^{p^\ast}} \leq C_{p,N} \|\nabla u\|_{L^p} \qquad \text{for every} \quad u \in \dot W^{1,p}(\R^N),\]
and hence on these spaces we have a Poincaré type inequality, which is exactly what we need for our estimates to work also in unbounded settings. This also tells us that the norm \eqref{norm-dot} is equivalent to
\[\|u\|_{L^{p^\ast}} + \|\nabla u\|_{L^p}.\]
With these considerations, $\dot W^{1,p}(\R^N)$ is a reflexive Banach space (see \cite[Theorem 12.2.3]{homog-book}) and, therefore, \autoref{Stamp} can be applied in this setting.

Given $\varphi \in W^{1-\frac1{p},p}(\partial \Omega)$,  \cite[Theorem 18.40]{Leoni} ensures that the trace opera\-tor is surjective on sets with Lipschitz continuous boundary, 
thus we can find $\Phi \in W^{1,p}(\Omega)$ whose trace on $\partial\Omega$ is $\varphi$. Then, given $u\in L^{p^\ast}(\R^N\setminus\overline\Omega)$ with $\nabla u\in L^p(\R^N\setminus\overline\Omega)$ such that $u=\varphi$ on $\partial\Omega$, 
we will consider the extension
\[\tilde u := u \, \chi_{\R^N \setminus \Omega} + \Phi \,\chi_\Omega \in \dot{W}^{1,p}(\R^N).\]
Hence we work with functions on the space
\[\dot{W}^{1,p}_{\Phi}(\R^N) = \big\{v\in \dot{W}^{1,p}(\R^N) \, |\, v = \Phi \text{ in }\Omega\big\}.\]

Now we consider the energy functional
$$\mathcal{F}_{\R^N \setminus \overline \Omega, \Phi}:L^{p^*}(\R^N)\rightarrow [0,+\infty]$$  defined by
\begin{equation} \label{defF_hom}
\mathcal{F}_{\R^N \setminus \overline \Omega,\Phi}(u):=\left\{\begin{array}{ll}
\displaystyle\int_{\R^N} F^p(\nabla u), & \hbox{if } u\in \dot W^{1,p}_{\Phi}(\R^N),   \smallskip \\
+\infty, & \hbox{otherwise},
\end{array}\right.
\end{equation}
and observe that it is a proper and convex functional. To characterize its subdifferential, we proceed exactly as before, that is, we introduce an auxilia\-ry operator $\mathcal A_\Phi \subset L^{p^\ast}(\R^N) \times  L^{(p^\ast)'}(\R^N)$ by doing the obvious modifications in Definition \ref{defA}, and repeat the arguments from the previous section for a sequence of approximating operators $\mathcal A_\Phi^{n,m}: \dot W^{1, p}_\Phi(\R^N) \rightarrow  \big(\dot W^{1, p}(\R^N)\big)'$ given again by the expression \eqref{def-aproxA}. The aforementioned properties of homogeneous Sobolev spaces allow us to repeat almost verbatim the steps in the previous section to conclude the following result:
\begin{theorem}[Characterization of the subdifferential on $\R^N \setminus \overline\Omega$] \label{Emin2}
	Let $1 < p < N$, $\Omega \subset \R^N$ an open bounded set with Lipschitz boundary and  $\Phi \in W^{1,p}(\Omega)$. If $u\in \dot{W}^{1,p}_\Phi(\R^N)$ 
and  $v\in L^{(p^\ast)'}(\R^N)$, then the following are equivalent:
	\begin{itemize}
		\item $v  \in \partial\mathcal{F}_{\R^N \setminus \overline \Omega,\Phi}(u)$.
		\item There exists $z\in L^\infty(\R^N\setminus \overline{\Omega};\R^N)$, with $F^\circ(z) \leq 1$ and $z \cdot \nabla u = F(\nabla u)$, such that $v=-{\rm div} (p F^{p-1}(\nabla u) z)$ in the weak sense, that is,
		$$\int_{\R^N \setminus \overline \Omega} vw = \int_{\R^N \setminus \overline \Omega} p F^{p-1}(\nabla u) z\cdot\nabla w \quad \hbox{ for every } w\in \dot W_{0}^{1,p}(\R^N),$$where $\dot W_{0}^{1,p}(\R^N)=\big\{v\in \dot{W}^{1,p}(\R^N) \, : \, v = 0 \text{ in }\Omega\big\}$.
	\end{itemize}
\end{theorem}

We notice that Theorem \ref{Emin2} gives a characterization for all minimizers of ${\rm Cap}_p^F(\overline\Omega)$ for bounded Lipschitz domains $\Omega$ since $${\rm Cap}_p^F(\overline\Omega)=\inf\left\{\mathcal F_{\R^N\setminus\overline\Omega,1}(u) : u\in L^{p^*}(\R^N)\right\}.$$

\section{A comparison principle for strictly convex norms and a barrier argument} \label{comp-barr}

\subsection{Results for strictly convex anisotropies}

 We start by obtaining a comparison principle for problem \eqref{EquationBounded} with a strictly convex even anisotropy $F$. It is well-known (cf.~\cite[Corollary 1.7.3]{Schneider}) that strict convexity of $F$ is closely related to differentiability of $F^\circ$. More precisely,
 \[F^\circ \in C^1(\R^N\setminus \{0\}) \,\, \text{if and only if} \,\, B_F(1)=\{\xi: F(\xi)<1\} \text{ is strictly convex}.\]
In the above conditions, it may happen that $F \not\in C^1(\R^N\setminus \{0\})$, as shown by the examples in \cite[Example A.1.19]{CaSin} or \cite[Section A]{CiGre}. In this spirit, the following result generalizes \cite[Lemma 2.3]{BiCi} by relaxing the regularity requirements for $F$ and for the functions $u_i$ involved.
\begin{lemma}\label{ComparisonPrincipleForApproximations}
	Let $\Omega\subset\R^N$  be a bounded domain with Lipschitz continuous boundary and $F$ a strictly convex norm in $\R^N$. If $u_1, u_2 \in W^{1,p}(\Omega)$ 
satisfy
	\[\left\{ \begin{array}{ll} {\rm div }\big(F^{p-1}(\nabla u_1) z_1\big) \geq {\rm div }\big(F^{p-1}(\nabla u_2) z_2\big) & \text{in }  \Omega \smallskip \\
	u_1 \leq u_2 & \text{on } \partial\Omega
	\end{array} \right.\]
	in the weak sense for some $z_i \in \partial F(\nabla u_i)$, $i =1,2$, then $u_1 \leq u_2$ in $\Omega$.	
\end{lemma}

\begin{proof}
	Setting $\Omega^+ = \Omega \cap \{u_1 > u_2\}$ and using $(u_1-u_2)^+$ as a test function, we have
	\[\int_{\Omega^+} \big(F^{p-1}(\nabla u_1) z_1 -F^{p-1}(\nabla u_2) z_2\big)\cdot \nabla(u_1 - u_2)\leq 0.\]	
	Notice that no boundary terms arise since $(u_1 - u_2)^+ = 0$ on $\partial\Omega^+$.
	
	On the other hand, as  $p F^{p-1}(\nabla u_i) z_i \in \partial F^p(\nabla u_i)$ and the latter is monotone,  we reach
	\[\big(F^{p-1}(\nabla u_1) z_1 -F^{p-1}(\nabla u_2) z_2\big)\cdot \big(\nabla u_1 - \nabla u_2\big)\geq 0.\]
	Therefore, the integrand on left hand side vanishes a.e. in $\Omega^+$. Accordingly, by means of \eqref{char-dF}, we get
	\begin{align*}
	0 & = F^p(\nabla u_1) + F^p(\nabla u_2) - F^{p-1}(\nabla u_1) z_1 \cdot \nabla u_2 - F^{p-1}(\nabla u_2) z_2 \cdot \nabla u_1
	\\ & \geq \big(F^{p-1}(\nabla u_1) - F^{p-1}(\nabla u_2)\big)  \big(F(\nabla u_1) - F(\nabla u_2)\big),
	\end{align*}	
	meaning that $F(\nabla u_1) = F(\nabla u_2)$ a.e. in $\Omega^+$. Hence the strict convexity of $F$ ensures that $\nabla u_1 = \nabla u_2$ a.e. in $\Omega^+$. As $u_1 = u_2$ on $\partial \Omega^+$, we conclude that $\Omega^+$ has null measure and thus the claim follows.
\end{proof}

\begin{remark}
	We observe that this result gives uniqueness of minimizers  of $\mathcal F_{\Omega,\varphi}$ (and solutions to \eqref{EquationBounded}) if $F$ is a strictly convex norm.
\end{remark}

We also need the following auxiliary result:
\begin{lemma} \label{F-strict}
Let $F$ be a strictly convex norm in $\R^N$, then for every $x \neq 0$
\[F(\nabla F^\circ(x)) = 1 \qquad \text{and}\qquad
\frac{x}{F^\circ(x)} \in \partial F(\nabla F^\circ(x)).\]
\end{lemma}

\begin{proof}
The proof of the first equality is exactly as in \cite[Proposition 1.3]{Xia}, since it only uses the differentiability of $F^\circ$, which follows from the strict convexity of $F$. For the second claim, from the 1-homogeneity of $F^\circ$ (recall \eqref{char-dF}), we can write
\[ x \cdot\nabla F^\circ(x) = F^\circ(x) =F(\nabla F^\circ(x)) F^\circ(x).\]
By \eqref{CS-ineq} this implies that $\nabla F^\circ(x)$ minimizes  $g(\xi):= F^\circ(x) F(\xi) - x \cdot \xi$. Accordingly,
\[0 \in \partial g(\nabla F^\circ(x)) = F^\circ(x) \partial F(\nabla F^\circ(x))-x,\]
meaning that $x \in F^\circ(x) \partial F(\nabla F^\circ(x))$, as desired.
\end{proof}

\subsection{Existence of minimizers within annular domains}

\begin{notation} \label{not-R}
	For any $R>0$ and $\Omega$ a bounded domain with Lipschitz-continuous bounda\-ry in $\R^N$ such that $\overline\Omega\subset \mathcal W_R$, we consider the annular region $\Omega_R:= \mathcal W_R \setminus \overline{\Omega}$ and the function $\psi$  defined on $\partial\Omega_R$ by $\psi=1$ on $\partial\Omega$ and $\psi=0$ on $\partial \mathcal W_R$. Now it makes sense to work with the functional $\mathcal{F}_{\!\!_R}:=\mathcal{F}_{\Omega_R,\psi}$ given by \eqref{defF}.
\end{notation}

In this setting, we aim to obtain a minimizer of  $\mathcal{F}_{\!\!_R}$ and then apply suita\-ble comparison results to take limits as $R\to \infty$. 	However, in order to exploit Lemma \ref{ComparisonPrincipleForApproximations}  in general, that is, for non-necessarily strictly convex norms, we will approximate the anisotropy $F$ by a sequence of strictly convex anisotropies $\{F_\varepsilon = F + \varepsilon \|\cdot \|\}_{\varepsilon >0}$, and consider the subsequent energy functional
$$\mathcal{F}_{\!\!_R}^\varepsilon:L^p(\Omega_R)\rightarrow [0,+\infty]$$  defined by
$$\mathcal{F}_{\!\!_R}^\varepsilon(u):=\left\{\begin{array}{ll}
\displaystyle\int_{\Omega_R} F_\varepsilon^p(\nabla u) & \hbox{if } u\in W^{1,p}_\psi(\Omega_R) \smallskip \\ 
+\infty & \hbox{otherwise}
\end{array}\right. ,$$
 where, by simplicity,  we  denote $(F_\varepsilon)^p$ by $F_\varepsilon^p$.

\begin{lemma} \label{seq-min}
Let $\{u_\varepsilon\}_{\varepsilon >0}$ be a sequence of minimizers of $\mathcal{F}_{\!\!_R}^\varepsilon$. Then there exists a subsequence that converges weakly in $W^{1,p}(\Omega_R)$ to $u_R \in W^{1,p}_\psi(\Omega_R)$. Moreover, $u_R$ is a minimizer of $\mathcal{F}_{\!\!_R}$.
\end{lemma}

\begin{proof}
Let $u_\varepsilon$ be a minimizer of $\mathcal{F}_{\!\!_R}^\varepsilon$ and let $w\in W^{1,p}_\psi(\Omega_R)$. By the uniform convergence of $\{F_\varepsilon\}_{\varepsilon>0}$ to $F$ as $\varepsilon\searrow0$ and \eqref{FequivEuclid}, there exist $k>0$ not depending on $\varepsilon$ such that  
$$(C+1)^p\|\nabla w\|_{L^p}^p\ge \mathcal{F}_{\!\!_R}^\varepsilon(w) \ge \mathcal{F}_{\!\!_R}^\varepsilon(u_\varepsilon)\ge k \Vert \nabla u_\varepsilon\Vert_{L^p}^p \quad \hbox{for all } \varepsilon \hbox{ small enough.}$$
Accordingly, by Poincar\'e inequality and up to a subsequence, we may suppose that $u_\varepsilon\rightharpoonup u_R\in W^{1,p}(\Omega_R)$ weakly in $W^{1,p}(\Omega_R)$ and strongly in $L^p(\Omega_R)$.

Now, since $ \mathcal{F}_{\!\!_R}$ is lower semicontinuous with respect to the weak convergence in $W^{1,p}(\Omega_R)$ (recall Proposition \ref{Flsc}), we can write
\begin{align*}
	\int_{\Omega_R} F^p(\nabla u_R) & \leq \liminf_{\varepsilon \to 0} 	\int_{\Omega_R} F^p(\nabla u_\varepsilon) \leq \limsup_{\varepsilon \to 0}  \int_{\Omega_R} \big(F(\nabla u_\varepsilon) + \varepsilon \|\nabla u_\varepsilon\|\big)^p \\ & \leq \limsup_{\varepsilon \to 0}  \int_{\Omega_R} \big(F(\nabla w) + \varepsilon \|\nabla w\|\big)^p = \int_{\Omega_R} F^p(\nabla w)
\end{align*}
for every $w \in W^{1,p}(\Omega_R)$. Therefore, we get that $u_R$ is a minimizer of $\mathcal{F}_{\!\!_R}$. By the continuity of  traces, $u_R = 1$ on $\partial \Omega$ and equals zero on $\partial \mathcal W_R$.
\end{proof}

\subsection{Construction of explicit barriers trapping a minimizer}

\begin{proposition} \label{comp-res}
	 Let $0<r_1<r_2 <R$ be such that $\mathcal W_{r_1}\subseteq \Omega\subseteq \mathcal W_{r_2}$, and let $u_R$ be a minimizer of $\mathcal{F}_{\!\!_R}$ which arises as a subsequential limit of minimizers of $\{\mathcal{F}_{\!\!_R}^\varepsilon\}_{\varepsilon>0}$ as $\varepsilon \searrow 0$. Then
	 \begin{equation}\label{boundsuR}\frac{\left(F^\circ\right)^{\frac{p-N}{p-1}}(x)-R^{\frac{p-N}{p-1}}}{r_1^{\frac{p-N}{p-1}}-R^{\frac{p-N}{p-1}}}\le u_R(x) \le \frac{\left(F^\circ\right)^{\frac{p-N}{p-1}}(x)-R^{\frac{p-N}{p-1}}}{r_2^{\frac{p-N}{p-1}}-R^{\frac{p-N}{p-1}}} \end{equation}
	 for a.e. $x\in \Omega_R$.
\end{proposition}

\begin{proof}
Given $r<R$, it is a routine computation to check that
$$v^\varepsilon_{r,R}(x)=\frac{\left(F_\varepsilon^\circ\right)^{\frac{p-N}{p-1}}(x)-R^{\frac{p-N}{p-1}}}{r^{\frac{p-N}{p-1}}-R^{\frac{p-N}{p-1}}}$$
satisfies
\[\left\{ \begin{array}{ll} {\rm div }\Big(F^{p-1}_\varepsilon(\nabla v^\varepsilon_{r,R}) \frac{x}{F_\varepsilon^\circ(x)}\Big) =0 & \text{in }  \mathcal W^\varepsilon_R\setminus\overline{\mathcal W^\varepsilon_r} \smallskip \\
v^\varepsilon_{r,R} = 1 & \text{on } \partial \mathcal W^\varepsilon_r \smallskip \\
v^\varepsilon_{r,R}=0 & \text{on } \partial \mathcal W^\varepsilon_R
\end{array} \right. .\]
Here $\mathcal W^\varepsilon_r$ is the Wulff shape of radius $r$ associated to the anisotropy $F_\varepsilon$.   Notice that the first equality follows after checking that $\partial F_\varepsilon(\nabla v^\varepsilon_{r,R}) = \partial F(\nabla F_\varepsilon^\circ)$ and using Lemma \ref{F-strict}. It can be easily checked, as in \cite[Lemma 5.3]{CCMN} that \begin{equation}\label{Wulff} \mathcal  W^\varepsilon_r = \mathcal W_r+B_{\|\cdot\|}(\varepsilon r).\end{equation}

On the other hand, as $u_\varepsilon\in W^{1,p}_\psi(\Omega_R)$ 
 is a minimizer for $\mathcal{F}_{\!\!_R}^\varepsilon$, then it holds
\[\left\{ \begin{array}{ll} {\rm div }\big( F^{p-1}_\varepsilon(\nabla u_\varepsilon) z\big) =0 & \text{in }  \Omega_R \smallskip \\
u_\varepsilon = 1 & \text{on } \partial \Omega \smallskip \\
u_\varepsilon=0 & \text{on } \partial \mathcal W_R
\end{array} \right., \qquad \text{with} \ z \in \partial F_\varepsilon(\nabla u_\varepsilon).\]

Let $r_2< R_1< R$ be such that $\mathcal W^\varepsilon_{r_1}\subseteq \Omega \subseteq \mathcal W^\varepsilon_{r_2}\subseteq \mathcal W^\varepsilon_{R_1}\subseteq \mathcal W_R$. Notice that $v^\varepsilon_{r_1, R_1} \leq u_\varepsilon$ on $\partial \Omega_R$, where $v^\varepsilon_{r_1, R_1}$ is extended by zero in $\mathcal W_R \setminus \mathcal W_{R_1}^\varepsilon$.   Therefore Lemma \ref{ComparisonPrincipleForApproximations} ensures that
\begin{equation} \label{barrier_up}
u_\varepsilon(x) \geq v^\varepsilon_{r_1, R_1}(x)= \frac{\left(F_\varepsilon^\circ\right)^{\frac{p-N}{p-1}}(x)-R_1^{\frac{p-N}{p-1}}}{r_1^{\frac{p-N}{p-1}}-R_1^{\frac{p-N}{p-1}}} \quad \text{for a.e. }\ x\in \Omega_R.
\end{equation}

Similarly, for $R_2 >R$ with $\mathcal W_R \subseteq \mathcal W^\varepsilon_{R_2}$, if we extend $v^\varepsilon_{r_2, R_2}$ by 1 in $\mathcal W_{r_2}^\varepsilon \setminus \overline \Omega$, we have that $u_\varepsilon\leq v^\varepsilon_{r_2, R_2}$ on $\partial \Omega_R$, and hence
\begin{equation} \label{barrier_down}
u_\varepsilon(x) \leq v^\varepsilon_{r_2, R_2}(x)= \frac{\left(F_\varepsilon^\circ\right)^{\frac{p-N}{p-1}}(x)-R_2^{\frac{p-N}{p-1}}}{r_2^{\frac{p-N}{p-1}}-R_2^{\frac{p-N}{p-1}}} \quad \text{for a.e. }\ x\in \Omega_R.
\end{equation}

Finally,  given $\varepsilon>0$, we set
\[R_1(\varepsilon)=\max \{\rho >0 \, | \, \mathcal W_\rho^\varepsilon\subseteq \mathcal W_R \} \ \text{ and } \ R_2(\varepsilon)=\min \{\rho >0 \, | \, \mathcal W_{R}\subseteq \mathcal W_{\rho}^\varepsilon\}.\]
 Then, by \eqref{Wulff}, we have
 \[R \leq \lim_{\varepsilon\searrow 0}R_1(\varepsilon)\leq\lim_{\varepsilon\searrow 0}R_2(\varepsilon)\leq R.\] So both $R_1$ and $R_2$ converge to $R$ and, taking limits in \eqref{barrier_up} and \eqref{barrier_down}, we conclude \eqref{boundsuR}, as desired.
\end{proof}

\section{Existence of minimizers with double side bounds} \label{sec-ult}



\begin{proof}[Proof of \autoref{umin-trapped}]
	
	Recall that for any $R>0$ one can actually construct a minimizer $u_R$, by means of Lemma \ref{seq-min}, for which  we have upper and lower barriers as stated in Proposition \ref{comp-res}. Now the aim is to show that these estimates pass to the limit as $R\to \infty$.

	Let $0<r_1<r_2$ such that $\mathcal W_{r_1}\subseteq \Omega\subseteq \mathcal W_{r_2}$. For each big enough $R$, let $u_R$ be a minimizer of $\mathcal{F}_{\!\!_R}$ which arises as a limit of minimizers of $\{\mathcal{F}_{\!\!_R}^\varepsilon\}_{\varepsilon>0}$ and extend $u_R$ by zero to $\R^N\setminus\overline{\Omega}$ (without renaming). Then by Lemma \ref{ComparisonPrincipleForApproximations} we have that $\{u_R\}$ is increasing in $R$ and, by Proposition \ref{comp-res}, we get
	\[\frac{\left(F^\circ\right)^{\frac{p-N}{p-1}}(x)-R^{\frac{p-N}{p-1}}}{r_1^{\frac{p-N}{p-1}}-R^{\frac{p-N}{p-1}}}\le u_R(x) \le \frac{\left(F^\circ\right)^{\frac{p-N}{p-1}}(x)-R^{\frac{p-N}{p-1}}}{r_2^{\frac{p-N}{p-1}}-R^{\frac{p-N}{p-1}}} \quad
	\text{a.e. in } \R^N\setminus \overline{\Omega}. \]
	Consequently, there exists a measurable function $u:\R^N\setminus \overline{\Omega}\to\R$ such that, up to a subsequence, $u_R(x)\nearrow u(x)$ as $R\to\infty$ for a.e. $x\in\R^N\setminus \overline{\Omega}$ and
	\begin{equation}\label{boundsu}\left(F^\circ\right)^{\frac{p-N}{p-1}}\bigg(\frac{x}{r_1}\bigg)\le u(x) \le \left(F^\circ\right)^{\frac{p-N}{p-1}}\bigg(\frac{x}{r_2}\bigg) \quad \hbox{a.e. in } \R^N\setminus\overline{\Omega}.
	\end{equation}
	In particular, $u\to0$ as $\|x\|\to\infty$. On the other hand, notice that the bounding functions in \eqref{boundsu} belong to $L^q(\R^N\setminus\overline\Omega)$ if $q > \frac{N(p-1)}{N-p}$, and this ensures by dominated convergence that $u_R\to u$ in $L^{p^\ast}(\R^N\setminus\overline{\Omega})$.
	
	Now, let $\varphi\in C_c^\infty (\R^N)$  with $\varphi=1$ on $\partial\Omega$, and consider  $R$ large enough so that $\mbox{supp}\, \varphi\subset \mathcal W_R$. Since $u_R$ is a minimizer of $\mathcal{F}_{\!\!_R}$, we have that
	$$c^p\Vert \nabla u_R \Vert^p_{L^p(\R^N\setminus\overline{\Omega})}\le \mathcal{F}_{\!\!_R}(u_R)\le \mathcal{F}_{\!\!_R}(\varphi)\le C^p\Vert \nabla\varphi\Vert^p_{L^p(\R^N\setminus\overline{\Omega})}.$$
	Therefore, $\{\Vert\nabla u_R \Vert_p\}_R$ is bounded thus $u\in W^{1,p}(\R^N\setminus\overline{\Omega})$ and $\nabla u_R\rightharpoonup \nabla u$ weakly in $L^p(\R^N\setminus\overline{\Omega})$; in particular, $u=1$ on $\partial\Omega$. Moreover, since $0\in\partial\mathcal{F}_{\!\!_R}(u_R)$, by Theorem \ref{Emin} we have that there exists $z_R\in\partial F(\nabla u_R)$ (which we extend by zero in $\R^N\setminus\mathcal W_R$) such that
	\begin{equation}\label{DivZR0} \int_{\R^N\setminus\overline{\Omega}} pF^{p-1}(\nabla u_R)z_R\cdot\nabla w=0
	\end{equation}
	for every $w\in W_0^{1,p}(\mathcal W_R\setminus\overline{\Omega})$. Now, by \eqref{subdif-bdd}, we have $\Vert z_R\Vert_\infty\le C$ so, up to a subsequence,
	\begin{equation}\label{ConvOfzRweakLpprime}
		pF^{p-1}(\nabla u_R)z_R\rightharpoonup \tilde{z} \ \ \hbox{in the weak topology of } L^{p'}(\R^N\setminus\overline{\Omega}).
	\end{equation}
	Taking limits as $R\to\infty$ in \eqref{DivZR0} we get
	$$ \int_{\R^N\setminus\overline{\Omega}} \tilde{z}\cdot\nabla w=0$$
	for every $w\in C_c^\infty(\R^N\setminus\overline{\Omega})$ (thus, by approximation, for any $w\in W_0^{1,p}(\R^N\setminus\overline{\Omega})$ that is,
	\begin{equation}\label{DivOfZRtilde}
		-{\rm div} (\tilde{z})=0 \ \ \hbox{ in the weak sense}.
	\end{equation}
	
	It remains to show that $\tilde{z}\in\partial F^p(\nabla u)$, which will be deduced 
	as usual from Lemma \ref{Minty-tec}. Therefore, it suffices to get \eqref{Minty2} for $z = \tilde z$ and $\alpha = 0$. With this goal, let $\xi\in\R^N$ and take $\omega:\R^N \setminus \overline\Omega\rightarrow \R$ defined by $\omega(\eta):=\xi\cdot\eta$ so that $\nabla \omega =\xi$. Then, letting $\phi\in C^\infty_c(\R^N \setminus \overline\Omega)$ such that $\phi\ge 0$, by the monotonicity of $\partial F^p$ and \eqref{DivZR0} we have, for every $R >0$ and any $z' \in \partial F^p(\xi)$,
	$$\begin{array}{l}\displaystyle 0\le\int_{\R^N \setminus \overline\Omega}\left(pF^{p-1}(\nabla u_R)z_R-z'\right)\cdot (\nabla u_R-\nabla \omega)\phi\\ [12pt]
		\displaystyle \phantom{0}= -\int_{{\rm supp} \nabla \phi}pF^{p-1}(\nabla u_R)z_R\cdot\nabla\phi(u_R-\omega)-\int_{\R^N \setminus \overline \Omega}z'\cdot (\nabla u_R-\nabla \omega)\phi.
	\end{array}$$
	Since by \eqref{boundsu} we have that $u_R \to u$ strongly in $L^p(K)$ for every compact set $K \subset \R^N \setminus \overline \Omega$, letting $R\to\infty$, by \eqref{ConvOfzRweakLpprime} and \eqref{DivOfZRtilde}, we get
	\begin{align*}\displaystyle 0 &\le \int_{\R^N \setminus \overline \Omega} -\tilde{z} \cdot \nabla \phi (u - \omega) -\int_{\R^N \setminus \overline \Omega}z'\cdot (\nabla u-\nabla \omega)\phi \\ & = \int_{\R^N \setminus \overline \Omega} \left(\tilde{z}-z'\right)\cdot(\nabla u-\nabla\omega)\phi.
	\end{align*}
	Finally, since $0\le\phi\in C^\infty_c(\R^N \setminus \overline\Omega)$ is arbitrary and $\nabla \omega = \xi$ we reach \eqref{Minty2} as desired.
\end{proof}

\section{Lipschitz regularity of minimizers for domains satisfying a uniform interior ball condition} \label{reg-sect}


If we further assume that $\Omega$ satisfies the condition stated in Definition \ref{UIBC}, by a careful application of our comparison results, we will show that a minimizer constructed as in Lemma \ref{seq-min} is bounded above and below by two Lipschitz solutions which coincide on the boundary, and accordingly we will deduce the Lipschitz regularity of the minimizer by application of the following result (see \cite[Corollary 4.2]{MT} for a more general version):
\begin{theorem} \label{MaTreu}
Let $f: \R^N \rightarrow \R$ be a convex function which is bounded below by an affine function and let
$$W^{1,p}_{u_1,u_2}(\Omega)=\{u\in W^{1,p}(\Omega) \, : \, u_1\le u\le u_2 \hbox{ a.e. and } u=u_1=u_2 \hbox{ on } \partial\Omega \},$$
where $u_1$ and $u_2$ are Lipschitz functions on $\Omega$ coinciding on $\partial\Omega$.
 Then the problem of minimizing $I(u) = \int_\Omega f(\nabla u)$ in $W^{1,p}_{u_1,u_2}(\Omega)$ admits at least one solution; moreover, at least one of the minimizers is Lipschitz with Lipschitz constant $L\le\max\{{\rm Lip}(u_1), {\rm Lip}(u_2)\}$.	
\end{theorem}

 With the conventions from Notation \ref{not-R}, we are now in a position to prove that any minimizer $w_R$ of the energy functional $\mathcal{F}_{\!\!_R}$  is Lipschitz continuous. Moreover,  any minimizer $u$ of $\mathcal{F}_{\R^N \setminus{\overline \Omega},1}$  is Lipschitz continuous.


\begin{proof}[Proof of Theorem \ref{lip-thm}]
	Let $R_0:=r+\frac{1}{2}{\rm dist }_{F^\circ}(\partial\Omega,\partial \mathcal W_R)$. For each $z\in\partial\Omega$ let $y_z\in\R^N$ such that $\mathcal W_r+y_z\subseteq \overline{\Omega}$ and $z\in\partial \left(\mathcal W_r+y_z\right)$. Given $z\in\partial\Omega$, define $\eta_z : \R^N\setminus \{y_z\} \rightarrow \R$ by
	$$\eta_z(x):=\frac{\left(F^\circ\right)^{\frac{p-N}{p-1}}(x-y_z)-R_0^{\frac{p-N}{p-1}}}{r^{\frac{p-N}{p-1}}-R_0^{\frac{p-N}{p-1}}}.$$
	Note that, $\eta_z(z)=1$, $\eta_z\le 1$ on $\partial \Omega$ and 
 $\eta_z\le 0$ on $\partial \mathcal W_R$. Moreover, as in the proof of Proposition \ref{comp-res}, one shows that $\eta_z(x)\le u_R(x)$ for a.e. $x\in \Omega_R= \mathcal W_R\setminus\overline{\Omega}$ and every $z\in\partial\Omega$. In addition, since $F^\circ(x-y_z)\ge r$ for every $x\in \Omega_R$ and
	$$\partial \eta_z(x) :=\frac{\frac{p-N}{p-1}\left(F^\circ\right)^{\frac{1-N}{p-1}}(x-y_z)\partial F^\circ(x-y_z)}{r^{\frac{p-N}{p-1}}-R_0^{\frac{p-N}{p-1}}},$$
	we have, for every $x\in \Omega_R$ and $w\in \partial \eta_z(x)$, that
	$$\Vert w\Vert\le  \frac{\frac{|N-p|}{p-1}r^{\frac{1-N}{p-1}}}{r^{\frac{p-N}{p-1}}-R_0^{\frac{p-N}{p-1}}} \,\Vert \tilde w\Vert \qquad \text{ where } \tilde w \in \partial F^\circ(x-y_z).$$
	Hence by \eqref{F0equivEuclid} and \eqref{subdif-bdd} we get the bound
	$$\Vert w\Vert\le  \frac{\frac{|N-p|}{p-1}r^{\frac{1-N}{p-1}}}{r^{\frac{p-N}{p-1}}-R_0^{\frac{p-N}{p-1}}} \frac{1}{c}=:L_1$$
	Therefore, $\eta_z$ is $L_1$-Lipschitz for every $z\in\partial \Omega$ (with $L_1$ independent of $R$).

	Let $I$ be a countable dense subset of $\partial\Omega$ and
	$$\eta(x):=\sup_{z\in I}\left\{\eta_z(x)\vee 0\right\}.$$
	Since the $\eta_z$ are $L_1$-Lipschitz we get that $\eta$ is also $L_1$-Lipschitz. Moreover, $\eta=1$ on $\partial\Omega$, $\eta=0$ on $\partial \mathcal W_R$ and $\eta\le u_R$  a.e. in $\Omega_R$.
	
	Therefore, if $\mathcal W_{r_1} \subseteq \Omega \subseteq \mathcal W_{r_2}\subseteq \mathcal W_R$, Proposition \ref{comp-res} leads to
	$$\eta(x)\le u_R(x)\le \frac{\left(F^\circ\right)^{\frac{p-N}{p-1}}(x)-R^{\frac{p-N}{p-1}}}{r_2^{\frac{p-N}{p-1}}-R^{\frac{p-N}{p-1}}}\wedge 1=: \eta'(x) \quad \hbox{for a.e. } x\in \Omega_R.$$
	Working as above we get that $\eta'$ is $L_2$-Lipschitz for some constant $L_2$ which decreases with $R$. In short, we have shown that $u_R \in W^{1,p}_{\eta, \eta'}(\Omega_R)$.

	Now, by \autoref{MaTreu}, there exists at least one minimizer of $\mathcal{F}_{\!\!_R}$ within $W^{1,p}_{\eta, \eta'}(\Omega_R)$ which is Lipschitz continuous with constant $L\le\max\{L_1,L_2\}$. If we denote this minimizer by $v_R$, notice that
	\[\mathcal{F}_{\!\!_R}(u_R) = \min_{W^{1,p}(\Omega_R)} \mathcal{F}_{\!\!_R} =  \min_{W^{1,p}_{\eta, \eta'}(\Omega_R)} \mathcal{F}_{\!\!_R} = \mathcal{F}_{\!\!_R}(v_R). \]
	
Next consider any minimizer $w_R$ of $\mathcal{F}_R$ (not necessarily constructed by a limiting procedure or trapped between $\eta$ and $\eta'$). Therefore,
	\begin{align}
		\mathcal{F}_{\!\!_R}\left(\frac{w_R+v_R}{2}\right)&=\int_{\Omega_R}F^p\left(\frac{\nabla(w_R+v_R)}{2}\right)\\
		&\le \int_{\Omega_R}\frac{F^p(\nabla w_R)+F^p(\nabla v_R)}{2}=\int_{\Omega_R}F^p\left(\nabla w_R\right) = \mathcal{F}_{\!\!_R}(w_R).
	\end{align}
	As the reverse inequality holds because $w_R\in\argmin \mathcal{F}_R$, we conclude
	$$F^p\left(\frac{\nabla(w_R+v_R)}{2}\right)=\frac{F^p(\nabla w_R)+F^p(\nabla v_R)}{2} \quad \hbox{a.e. in } \Omega_R.$$
	Now, if $F(\nabla w_R)\neq F(\nabla v_R)$, by the strict convexity of $t\rightarrow |t|^p$ and the convexity of $F$, we have that (a.e. in $\Omega_R$)
	\begin{align}
		\left(\frac{F(\nabla w_R)+F(\nabla v_R)}{2}\right)^p &< \frac{F^p(\nabla w_R)+F^p(\nabla v_R)}{2}=F^p\left(\frac{\nabla(w_R+v_R)}{2}\right)\\
		&\le \left(\frac{F(\nabla w_R)+F(\nabla v_R)}{2}\right)^p,
	\end{align}
	which is a contradiction. Therefore, $F(\nabla w_R)= F(\nabla v_R)$ a.e. in $\Omega_R$ thus $w_R$ is $\frac{C}{c}L$-Lipschitz (as $v_R$ is $L$-Lipschitz), where $c$ and $C$ come from \eqref{FequivEuclid}.
	
	Now, since the $u_R$ are uniformly Lipschitz  and $u_R\nearrow u$ by the proof of  \autoref{umin-trapped}, we get that $u$ is Lipschitz continuous (with Lipschitz constant $\frac{C}{c}L$). Now take $w$ an arbitrary minimizer of $\mathcal{F}_{\R^N \setminus{\overline \Omega},1}$ and redo the previous argument with $u$ and $w$ playing the role of $v_R$ and $w_R$, respectively. Then we conclude that $w$ is also Lipschitz continuous with constant $\frac{C^2}{c^2}L$.
	%
\end{proof}

\appendix

\section{Non-uniqueness of minimizers in the case $p = 1$} \label{nonU}

\renewcommand{\theequation}{A.\arabic{equation}}

In this appendix, we show an example of non-uniqueness for relative $1$-capacitary functions; i.e., minimizers of the relative $1$-capacity with respect to a ball $\mathcal W_R$. We take $N=2$ and $F(\xi)=\|\xi\|_1=|\xi_1|+|\xi_2|$ as the anisotropy. Let $\Omega:=B_1$ be the (Euclidean) unit ball in $\R^2$ centered at the origin and set $R=2$, thus $\mathcal W_2=\{\xi\in \R^2: \|\xi\|_\infty=\max\{|\xi_1|,|\xi_2|\}\leq 2\}$. Consider
$${\rm Cap}_1^{\|\cdot\|_1}(\overline{B_1};\mathcal W_2)=\inf\left\{\int_{\mathcal W_2}\|\nabla u\|_1\,dx \, : \, u \in C_0^\infty(\mathcal W_2), \, u\geq 1 {\rm \ in \ }B_1\right\}.$$
It turns out that, in general, due to the lack of compactness of minimizing sequences in $W^{1,1}(\mathcal W_R)$, minimizers of linear growth functionals with respect to the gradient, as that in ${\rm Cap}_1^{\|\cdot\|_1}$, can be functions that do not belong to the Sobolev space $W^{1,1}(\mathcal W_R)$. Therefore, one needs to consider the relaxed functional to $G:L^1(\mathcal W_2\setminus\overline B_1)\to [0,+\infty]$, $$G(u):=\left\{\begin{array}
  {cc} \displaystyle\int_{\mathcal W_2\setminus\overline{B_1}}\|\nabla u\|_1 & {\rm if \ } u\in W^{1,1}(\mathcal W_2\setminus\overline{B_1}): \left\{ \begin{array}{l} u=1 {\rm \ in \ }\partial B_1, \\ u=0 {\rm \ in \ }\partial \mathcal W_2 \end{array} \right. \smallskip \\ +\infty & {\rm otherwise}
\end{array}\right.,$$
which is (see \cite[Theorem 4]{M}) $$\mathcal G(u):=|Du|_1({\mathcal W_2\setminus\overline{B_1}})+\int_{\partial B_1}\left\|\nu^{B_1}\right\|_1|u-1|\,d\mathcal H^1+\int_{\partial\mathcal W_2}\left\|\nu^{\mathcal W_2}\right\|_1|u|\,d\mathcal H^1,$$
where $\mathcal H^1$ denotes the 1-dimensional Hausdorff measure in $\R^2$.
 Moreover, $\nu^U$ represents the unit exterior normal to the open set $U\subset\R^2$; and the term $|Du|_1(U)$ means the anisotropic total variation measure of the open set $U\subset\R^2$, defined as (see \cite[Definition 3.1]{AmarBellettini}) $$|Du|_1(U):=\sup\left\{\int_U \ u \, {\rm div}z \, : \, z\in X_1(U)\right\}, $$ with $$X_1(U):=\left\{z\in L^\infty(U;\R^2)\, :\, {\rm supp} \, z\subset U,\, {\rm div}z\in L^2(U), \,  \|z\|_\infty\leq 1 \, {\rm a.e. \, in \, } U\right\}.$$
We note that $|Du|_1(U)$ is finite if, and only if, $u$ is a bounded variation function in $U$. Similarly, $E\subset\R^2$ is a set of finite perimeter if, and only if, $${\rm Per}_1(E):=|D\chi_E|_1(\R^2)<+\infty,$$
where $\chi_E$ is the characteristic function of the measurable set $E$; i.e., $$\chi_E(x):=\left\{\begin{array}
  {cc} 1 & {\rm if \ } x\in E \\ 0 & {\rm if \ } x\notin E
\end{array}\right..$$
For a comprehensive treatment of bounded variation functions and sets of finite perimeter, we refer to \cite{ambrosio2000fbv}.

The characterization of the subdifferential of the energy functional $\mathcal G$ is included in \cite[Theorem 9]{M}. We will show that $0\in \partial \mathcal G(u)$ in the case that $u=\chi_E$ with $E\subseteq\mathcal W_1\setminus\overline{B_1}$ of finite perimeter such that $\partial B_1\subset\partial E$.

In this setting, the characterization of the subdifferential is much simpler than for a generic bounded variation function and we obtain that $0\in\partial \mathcal G(\chi_E)$ if, and only if, there exists $z\in L^\infty(\mathcal W_2\setminus\overline{B_1})$ such that $\|z\|_\infty\leq 1$, a.e. in $\mathcal W_2\setminus\overline{B_1}$, ${\rm div } z=0$ in the distributional sense and \begin{equation}\label{app:subdiff}\mathcal G(\chi_E)=|D\chi_E|_1(\mathcal W_2\setminus \overline{B_1})={\rm Per}_1(E)-{\rm Per}_1(B_1)={\rm Per}_1(B_1).\end{equation}
We observe that the first two equalities hold because of the particular assump\-tion on $E$.

We define $$z(x,y):=\left\{\begin{array}{cc} \displaystyle -\left(\frac{x}{|x|}, \frac{y}{|y|}\right) & {\rm if \ } (x,y)\in\mathcal W_1\setminus\overline{B_1} \medskip \\ \displaystyle -\frac{(x,y)}{\|(x,y)\|_\infty^2} & {\rm if \ }(x,y)\in \mathcal W_2\setminus \mathcal W_1\end{array}\right..$$

It is easy to show that the vector field $z$ satisfies $\|z\|_\infty\leq 1$ in $\mathcal W_2\setminus\overline{B_1}$ and ${\rm div} z=0$ in the distributional sense.

On the other hand, a routine computation ensures that for \eqref{app:subdiff} to hold, we just need that \begin{equation}\label{app:cond_bdry}\nu^E(x,y)\cdot z(x,y)=-1\quad  \mathcal H^1{\rm -a.e. \ in  \ }\partial E\setminus \partial B_1.\end{equation}

Since there are infinite sets of finite perimeter $E\subset \mathcal W_1\setminus\overline{B_1}$ satisfying \eqref{app:cond_bdry}, we conclude that there are infinitely many different $1$-capacitary functions. For instance,
$E_1:=\mathcal W_1\cap( B_{\|\cdot\|_1}(\sqrt{2})\setminus \overline{B_1})$ and $E_2:=\mathcal W_1\setminus\overline{B_1}$ are two different examples. We finally note that it can be also proved, though the proof requires the general characterization of the subdifferential, that $u=\chi_\emptyset=0$ is a minimizer.

\end{document}